\documentclass[a4paper,11 pt]{article}
\usepackage[utf8]{inputenc}
\usepackage[T1]{fontenc}
\usepackage{amsfonts}
\usepackage{amssymb}
\usepackage{amsmath}
\usepackage{amsthm}
\usepackage{multicol}
\usepackage{euscript}
\usepackage{mathtools}
\usepackage[round]{natbib}
\usepackage{hyperref}
\usepackage[normalem]{ulem}
\mathtoolsset{
  showonlyrefs,
  mathic 
}
\hypersetup{
  colorlinks=true,
  breaklinks=true,
  linkcolor=red,
  anchorcolor=black,
  citecolor=green,
  filecolor=blue,
  menucolor=red,
  urlcolor=red,
}
\DeclareMathOperator*{\argmin}{arg\,min}

\theoremstyle{plain}
\newtheorem{theorem}{Theorem}
\newtheorem{proposition}{Proposition}

\newtheorem{lemma}{Lemma}
\newtheorem{definition}{Definition}

\theoremstyle{remark}
\newtheorem{remark}{Remark}

\let\fontRing\mathbb

\newcommand*{\R}{\fontRing{R}}

\newcommand*{\N}{\fontRing{N}} 

\let\fontProba\mathbf

\newcommand*{\p}{\fontProba{P}}
\newcommand*{\e}{\fontProba{E}}
\newcommand*{\1}{\fontProba{I}}

\newcommand*{\D}{{\,\mathrm{d}}}

\let\eps\varepsilon

\let\hat\widehat

\makeatletter
\def\calFactory#1{%
   \expandafter\def\csname c#1\endcsname{\mathcal{#1}}}
\def\frakFactory#1{%
   \expandafter\def\csname k#1\endcsname{\mathfrak{#1}}}
\def\bbFactory#1{%
   \expandafter\def\csname k#1\endcsname{\mathbb{#1}}}
\newcounter{ctr}
\loop
  \stepcounter{ctr}
  \edef\X{\@Alph\c@ctr}%
  \edef\x{\@alph\c@ctr}%
  \expandafter\calFactory\X
  \expandafter\frakFactory\X
  \expandafter\frakFactory\x
  \expandafter\bbFactory\X
\ifnum\thectr<26
\repeat
\makeatother

\newcommand{\dom}{\Delta}
\newcommand{\domain}[1]{\dom_{{#1}}}
\newcommand{\sobolev}{\mathbb{S}}
\newcommand{\lebesgue}{\mathbb{L}}

\newcommand{\Hn}{\mathcal H_n}

\newcommand{\const}[1]{\mathfrak{K}_{#1}}

\newcommand{\ms}{{\lfloor s\rfloor}}

\newcommand{\ellp}{{\ell'}}

\newcommand{\hell}{h(\ell)}

\newcommand{\h}{h}

\usepackage{xcolor}

\newcommand{\pfn}{\p_{\!\!f}^n}
\newcommand{\efn}{\e_f^n}
\newcommand{\Rnpq}{R_n^{(p,q)}}
\begin{document}

\title{Adaptive Density Estimation\\on Bounded Domains}
\author{
    Karine Bertin
    \thanks{CIMFAV, Universidad de
  Valpara\'{i}so, General Cruz 222,
  Valpara\'{i}so, Chile}
    \and
    Salima El Kolei
    \thanks{ENSAI, UBL, Campus de Ker Lann, Rue Blaise Pascal - BP 37203
35172 Bruz cedex}
    \and
    Nicolas Klutchnikoff
    \thanks{Univ Rennes, CNRS, IRMAR (Institut de Recherche Mathématique
de Rennes) - UMR 6625, F-35000 Rennes, France}
    }
\date{\today}
\maketitle

\begin{abstract}
    We study the estimation, in $\mathbb{L}_p$-norm, of density functions defined on $[0, 1]^d$. We construct a new family of kernel density estimators that do not suffer from the so-called boundary bias problem and we propose a data-driven procedure based on the Goldenshluger and Lepski approach that jointly selects a kernel and a bandwidth. We derive two estimators that satisfy oracle-type inequalities. They are also proved to be adaptive over a scale of anisotropic or isotropic Sobolev-Slobodetskii classes (which are particular cases of Besov or Sobolev classical classes). The main interest of the isotropic procedure is to obtain adaptive results without any restriction on the smoothness parameter.
{
    \footnotesize
    \par\medskip
    \noindent \textbf{Keywords.}  Multivariate kernel density estimation,  Bounded data, Boundary bias, Adaptive estimation, Oracle inequality, Sobolev-Slobodetskii classes.   \par\smallskip
    \noindent \textbf{AMS Subject Classification.} 62G05, 62G20.
}
\end{abstract}

\begin{abstract}
	Nous étudions l'estimation, en norme $\mathbb{L}_p$, d'une densité de probabilté définie sur $[0,1]^d$. Nous construisons une nouvelle famille d'estimateurs à noyaux qui ne sont pas biaisés au bord du domaine de définition et nous proposons une procédure de sélection simultanée d'un noyau et d'une fenêtre de lissage en adaptant la méthode développée par Goldenshluger et Lepski. Deux estimateurs différents, déduits de cette procédure générale, sont proposés et des inégalités oracles sont établies pour chacun d'eux. Ces inégalités permettent de prouver que les-dits estimateurs sont adapatatifs par rapport à des familles de classes de Sobolev-Slobodetskii anisotropes ou isotropes. Dans cette dernière situation aucune borne supérieure sur le paramètre de régularité n'est imposée.
\end{abstract}

\section{Introduction}

In this paper we study the classical problem of the estimation of a density function $f:\domain{d}\to\R$ where $\domain{d}=[0,1]^d$. We observe $n$ independent and identically distributed random variables $X_1,\ldots,X_n$ with density $f$. In this context, an estimator is a measurable map $\tilde f:\domain{d}^n\to \lebesgue_p(\domain{d})$ where $p\geq1$ is a fixed parameter. The accuracy of $\tilde f$ is measured using the risk:
\begin{equation}
  \Rnpq(\tilde f,f) = \left(\efn\|{\tilde f-f}\|_p^q\right)^{1/q},
\end{equation}
where  $q$ is also a fixed parameter greater than or equal to $1$ and $\efn$ denotes the expectation with respect to the probability measure $\pfn$ of the observations. Moreover the $\lebesgue_p$-norm of a function $g:\domain{d}\to\R$ is defined by
\begin{equation}
  \|g\|_p=\left(\int_{\domain{d}} |g(t)|^p dt\right)^{1/p}.
\end{equation}

We are interested in finding data-driven procedures that achieve the minimax rate of convergence over Sobolev-type functional classes that map $\Delta_d$ onto $\R$.
The density estimation problem is widely studied and we refer the reader to~\citet{devroye1985nonparametric} and~\citet{silverman1986density} for a broadly picture of this domain of statistics. One of the most popular ways to estimate a density function is to use kernel density estimates introduced by \citet{rosenblatt1956remarks} and \citet{parzen1962estimation}. Given a kernel $K$ (that is a function $K:\R^d\to\R$ such that $\int_{\R^d} K(x) \D x=1$) and a  bandwidth vector $h=(h_1,\ldots,h_d$), such an estimator writes:
\begin{equation}\label{eq:conv_kernel}
  \hat f_h(t) = \frac{1}{nV_h}\sum_{j=1}^n K\left(\frac{t-X_j}{h}\right), \qquad t\in\domain{d}
\end{equation}
where $V_h=\prod_{i=1}^d h_i$ and $u/v$ stands for the coordinate-wise division of the vectors $u$ and~$v$.

It is commonly admitted that \emph{bandwidth selection} is the main point to estimate accurately the density function $f$ and a lot of popular selection procedures are proposed in the literature. Among others let us point out the cross validation \citep[see][]{rudemo1982empirical,bowman1984alternative,chiu1991bandwidth} as well as the procedure developed by  Goldenshluger and Lepski in a series of papers in the last few years \citep[see][for instance]{MR2543590,GL2011,MR3230001} and fruitfully applied in various contexts.

Dealing with bounded data, the so-called \emph{boundary bias problem} has also to be taken into account. Indeed, classical kernels suffer from a severe bias term when the underlying density function does not vanish near the boundary of their support. To overcome this drawback, several procedures have been developed: \citet{schuster1985incorporating}, \citet{silverman1986density} and \citet{cline1991kernel} studied the reflection of the data near the boundary as well as \citet{marron1994transformations} who proposed a previous transformation of the data. \Citet{Muller1991Boundary},  \citet{LejeuneSarda1992Smooth}, \citet{jones1993simple}, \citet{muller1999multivariate} and \citet{botev2010kernel} proposed to construct kernels which take into account the shape of the support of the density. In the same spirit, \citet{MR1718494} studied a new class of kernels constructed using a reparametrization of the family of Beta distributions.
For these methods, practical choices of bandwidth or cross-validation selection procedures have generally been proposed. Nevertheless few papers study the theoretical properties of \emph{bandwidth selection} procedures in this context. Among others, we point out \citet{Bouezmarni2010Nonparametric}---who study the behavior of Beta kernels with a cross validation selection procedure in a multivariate setting in the specific case of a twice differentiable density. \Citet{bertin2014adaptive} study a selection rule based on the Lepski's method \citep[see][]{lepski1991asymptotically} in conjunction with Beta kernels in a univariate setting and prove that the associated estimator is adaptive over H\"older classes of smoothness smaller than or equal to two.
In this paper, we aim at constructing estimation procedures that address both problems (\emph{boundary bias} and \emph{bandwidth selection}) simultaneously and with optimal adaptive properties in $\mathbb{L}_p$ norm ($p\ge 1$) over a large scale of function classes. To tackle the boundary bias problem, we construct a family of kernel estimators based on new asymmetric kernels whose shape adapts to the position of the estimation point in $\domain{d}$. We propose two different data-driven procedures based on the Goldenshluger and Lepski approach that satisfy oracle-type inequalities (see Theorems~\ref{thm:oracle} and~\ref{thm:oracle2}). The first procedure, based on a fixed kernel, consists in selecting a bandwidth vector. Theorem~\ref{thm:adaptive-anisotropic} proves that the resulting estimator is adaptive over anisotropic Sobolev-Slobodetskii classes with smoothness parameters $(s_1,\ldots,s_d)\in(0,\infty)^d$ smaller than the order of the kernel and with the optimal rate $n^{-\overline{s}/(2\overline{s}+1)}$ with $\overline{s}=\left(\sum_{i=1}^d 1/s_i\right)^{-1}$. The second procedure jointly selects a kernel (and its order) and a univariate bandwidth. Such selection procedures have been used only in the context of exact asymptotics in pointwise and sup-norm risks, and for very restrictive function classes. Theorem~\ref{thm:adaptive-isotropic} states that this procedure is adaptive over isotropic Sobolev-Slobodetskii classes without any restriction on the smoothness parameter $s>0$ and achieves the optimal rate $n^{-s/(2s+d)}$.
These function classes are quite large and correspond to a special case of usual Besov classes \citep[see][]{MR1328645}. Note also that the same results can be obtained over anisotropic H\"older classes with the same rates of convergence.
Such adaptive results without restrictions on the smoothness of the function to be estimated and with the optimal rates $n^{-2s/(2s+d)}$ or $n^{-\overline{s}/(2\overline{s}+1)}$ have been established only for ellipsoid function classes as in~\cite{asin2016adaptive}, among others.
For bounded data, we also mention~\cite{rousseau2010rates} or~\cite{autin2010thresholding} that construct adaptive estimators based on Bayesian mixtures of Beta and wavelets respectively but with an extra logarithmic term factor in the rate of convergence.
Additionally note also that Beta kernel density estimators are minimax only for small smoothness \citep[see][]{bertin2011minimax} and consequently neither allow us to obtain such adaptive results.

The rest of the paper is organized as follows. In Section~\ref{sec:on_the_boundary_bias_problem}, we detail the effect of the boundary bias and we propose a new family of estimators that do not suffer from this drawback. We construct in Section~\ref{sec:stat} our two main statistical procedures. The main results of the paper are stated in Section~\ref{sec:results} whereas their proofs are postponed to Section~\ref{sec:proofs1}.

\section{On the boundary bias problem} 
\label{sec:on_the_boundary_bias_problem}

\subsection{Weakness of convolution kernel estimators}

In this section we focus on the so-called boundary bias problem that arises when classical convolution kernels are used. To illustrate our point and for the sake of simplicity we assume that $d=1$ and $p=q\geq 2$.
In what follows we consider the estimators defined in~\eqref{eq:conv_kernel}:
\[
  \hat{f}_h(t) = \frac{1}{nh} \sum_{j=1}^n K\left(\frac{t-X_j}{h}\right),
  \quad
  t\in\Delta_1
\]
where $0<h<1$ is a bandwidth and the kernel $K:\R\to\R$ is such that:
\[
  \operatorname{Supp}(K)\subseteq[-1,1]
  \quad\text{and}\quad
  \int_{-1}^0 K(u) du = 1-\gamma,
\]
with $0<\gamma<1$. In this context, the following lemma---which is straightforward---proves that these estimators suffer from an asymptotic pointwise bias at the endpoint $0$ as soon as $f(0)\neq 0$.
\begin{lemma}
  Assume that $f$ is continuous at $0$. Then, we have $\e_f^n\hat{f}_h(0) \xrightarrow[h\to0]{} (1-\gamma) f(0)$.
\end{lemma}
However this problem is not specific to the endpoint and generalizes to a whole neighborhood of this point. The simplest situation that allows one to understand this phenomenon is to consider the estimation of the function $f_0=\1_{(0,1)}$ where, here and after, $\1_{(a,b)}$ stands for the indicator function of the interval $(a,b)$.
In this case, under a more restrictive assumption on the kernel, the integrated bias can be bounded from below by $h^{1/p}$ up to a multiplicative factor. More precisely we can state the following result:
\begin{proposition}\label{thm:LB-bias}
    Assume that $K$ is a kernel such that $\operatorname{Supp}(K)\subseteq[-1,1]$ and assume that there exist $0<\delta<1$ and $0<\gamma<1$ that satisfy
    \[
      \int_{-1}^a K(u) du \leq 1-\gamma
      \quad\text{for any $0\leq a\leq\delta$.}
    \]
    Then, for any $0<h<(1+\delta)^{-1}$ we have:
    \[
      \|\e_{f_0}^n(\hat{f}_h)-f_0\|_p
      \geq
      \left(\delta^{1/p}\gamma\right)\, h^{1/p}.
    \]
\end{proposition}
As a consequence of this proposition we can state the following lower bound on the rate of convergence of the classical convolution kernel estimators over a very large family of functional classes.
\begin{proposition}\label{thm:LB-risk}
Let $p\ge 2$. Let $\Sigma$ be a functional class such that $f_0\in\Sigma$. Assume that $K\in\mathbb{L}_2([-1,1])$. Under the assumptions of Proposition~\ref{thm:LB-bias}, we have:
  \[
    \liminf_{n\to0}\inf_{h\in(0,1/4)} n^{1/(2+p)} \sup_{f\in\Sigma} R_n^{(p,p)}(\hat{f}_h,f) >0
  \]
\end{proposition}

Now, let us comment these two results. First, remark that the assumptions made on the kernel are not very restrictive since any continuous symmetric kernel $K$ such that $K(0)>0$ can be considered. Next, in view of Theorem~\ref{thm:adaptive-isotropic} stated below, Proposition~\ref{thm:LB-risk} proves that the convolution kernel estimators are not optimal. In particular, they do not achieve the minimax rate of convergence over usual H\"older classes with smoothness parameter $s>1/p$ (see Definition~\ref{def:ss-iso} as well as Remark~\ref{rem:main-results} for more details). This result is mainly explained by Proposition~\ref{thm:LB-bias} since, in this situation, the integrated bias term is greater than $h^{1/p}$ which is larger (in order) than the expected term $h^s$ (see Proposition~\ref{thm:UB-bias} below).

\subsection{Boundary kernel estimators}\label{ssec:1}

The main drawback of classical convolution kernels can be explained as follows: they \emph{look outside} the support of the function to be estimated.
As a consequence, $f_0$ is seen as a discontinuous function that maps $\R$ to $\R$. This leads to a severe bias and explains why ``boundary kernels'' found in the literature have all their mass \emph{inside} the support of the target function.  Indeed, in this situation $f_0$ is seen as a function that maps $\Delta_1$ to $\R$ which is a very smooth function. This allows the bias term to be small (see Proposition~\ref{thm:UB-bias} below)

In last decades, several papers proposed different constructions of kernels that can take into account the boundary problem. Let us point out that, among others, \citet{muller1999multivariate} and \citet{MR1718494} constructed specific kernels whose shape adapts to the localization of the estimating point in a continuous way. Even if this continuously deforming seems to be an attractive property there are still some drawbacks to using such approaches. On the one hand, the beta kernels cannot be used to estimate smooth functions~\citep[see][]{bertin2011minimax}. On the other hand, the kernels proposed by \citet{muller1999multivariate} are solutions of a continuous least square problem for each estimating point. In practice the kernels are computed using discretizations of the variational problems. This can be computationally intensive. Moreover, to our knowledge, there are no theoretical guarantees regarding bandwidth selection procedures in this context.

In this paper, we propose a simple and tractable way to construct boundary kernels that intends to solve the aforementioned problems. The main advantage of our construction lies in the fact that the resulting estimators are easy to compute and that the mathematical analysis of the adaptive procedure is made possible even in the anisotropic case.

To construct our kernels we first define the following set of univariate bounded kernels whose support is included into $\Delta_1$:
\begin{equation}
  \cW = \left\{W:\R\to\R : \sup_{u\in\domain1}|W(u)|<+\infty, \quad W(u)=0 \text{ for $u\notin \domain1$,} \quad \int_{\domain1} W(u) \D u=1\right\}.
\end{equation}
In the following, we will say that $W\in\cW$ is a kernel of order $m$ if
  \begin{equation}\label{eq:moment-condition}
    \int_{\domain1} W(u) u^r \D u = 0, \qquad r=1,2,\ldots,m.
  \end{equation}
Then, for any bandwidth $h\in(0,1/2)^d$ and any sequence of kernels $W=(W_1,\ldots,W_d)\in\cW^d$, we define the following density estimator:
\begin{equation}\label{def:fchapeau}
  \tilde f_{W,h}(t) = \frac1n\sum_{j=1}^n \mathcal K_{W,h}(t,X_j), \quad t\in\Delta_d
\end{equation}
where, for $t\in\Delta_d$ the ``boundary'' kernel $\mathcal{K}_{W,h}(t,\cdot)$ is defined by:
\begin{equation}
    \mathcal K_{W,h}(t,x) = \prod_{i=1}^d\left(\frac1h_i
    W_i\left(\sigma(t_i)\frac{t_i-x_i}{h_i}\right)
    \right)
\end{equation}
for any $x\in\Delta_d$. Here $\sigma(\cdot)=2\1_{(1/2,1)}(\cdot)-1$.

\begin{remark}
  Note that, along each coordinate, the kernel $W_i$ is simply flipped according to the position of $t_i$ with respect to the closest boundary. Similar constructions can be found in the literature. For example \citet{MR1735786} and \citet{MR2062980} proposed to decompose $\Delta_1$ into three different pieces --- that depend on the bandwidth $h$ --- as follows: $\Delta_1=(0,h)\cup[h,1-h]\cup(1-h,1)$. Specific kernels are used for the boundaries while classical kernels are used on $[h,1-h]$. However, to our best knowledge, similar constructions in a multivariate framework do not allow to obtain adaptive results in the anisotropic case.
  \end{remark}

\subsection{Bias over some functional classes}

In this paper we focus on minimax rates of convergence over Sobolev-Slobodetskii classes. We recall their definitions in Definitions~\ref{def1} and~\ref{def:ss-iso} \citep[see also][]{simon90,MR1120995,MR1328645}.

In the following, for \(f:\Delta_d\to\R\) and any \(i=1,\dotsc,d\) and $k\in\N$, we denote by \(D_i^k f\) the $k$th-order partial derivative of f with respect to the variable~$x_i$. For any $\alpha\in\N^d$, we denote by \(D^\alpha f\) the mixed partial derivatives
    \begin{equation}
        D^\alpha f = \frac{\partial^{|\alpha|} f}{\partial x_1^{\alpha_1}\cdots\partial x_d^{\alpha_d}}
    \end{equation}
where \(|\alpha|=\alpha_1+\dotsc+\alpha_d\). Finally, for any positive number $u$, we denote by $\lfloor u\rfloor$ the largest integer strictly smaller than $u$.

\begin{definition}\label{def1}
  Set $s=(s_1,\ldots, s_d)\in(0,+\infty)^d$  and $L>0$.
  A function $f:\domain{d}\to\R$, belongs to the anisotropic Sobolev-Slobodetskii ball $\sobolev_{p}(s,L)$ if:
  \begin{itemize}
    \item $f$ belongs to $\lebesgue_p(\domain{d})$.
    \item For any \(i=1,\dotsc,d\), $D_i^{\lfloor s_i\rfloor}f$ exists and belongs to $\lebesgue_p(\domain{d})$.
    \item The following property holds:
    \begin{equation}
      \sum_{i=1}^d I_i(D_i^{\lfloor s_i\rfloor}f) \leq L,
    \end{equation}
    where
    \begin{equation}
      I_i(g) = \left(\int_{\domain{d}}\int_{\domain{1}} \frac{\bigl|g(x)-g(x_1,\ldots,x_{i-1},\xi,x_{i+1},\ldots, x_d)\bigr|^p}{|x_i-\xi|^{1+p(s_i-\lfloor s_i\rfloor)}}dx d\xi\right)^{1/p}.
    \end{equation}
  \end{itemize}
\end{definition}

\begin{definition}\label{def:ss-iso}
Set $s>0$ and $L>0$.
A function $f:\domain{d}\to\R$, belongs to the isotropic Sobolev-Slobodetskii ball $\tilde\sobolev_{s,p}(L)$ if the following properties hold:
\begin{itemize}
    \item for any $\alpha\in\N^d$, such that $|\alpha|\leq\ms$, the mixed partial derivatives \(D^\alpha f \) exist and belong to $\lebesgue_p(\domain{d})$.
    \item the Gagliardo semi-norm $|f|_{s,p}$ is bounded by $L$ where
        \begin{equation}
        |f|_{s,p}=\left(\sum_{|\alpha|=\ms}\int_{\domain{d}^2} \frac{\left|{D^\alpha f(y) - D^\alpha f(x)}\right|^p}{\|{y-x}\|_2^{d+p(s-\ms)}} \D x\D y\right)^{1/p},
    \end{equation}
    where $\|\cdot\|_2$ denotes the euclidean norm of $\R^d$.
\end{itemize}
\end{definition}

These classes include several classical classes of functions. Indeed, in the isotropic case, when \(s>0\) is not an integer, then \(\tilde{\sobolev}_{s,p}(L)\) corresponds to the usual Besov ball $B_{p,p}^s(L)$ \citep[see][]{MR1328645}. Note that both definitions are the same when \(d=1\).

The following proposition illustrates the good properties in terms of bias of our boundary kernel estimators. It can be obtained following Propositions~\ref{prop:bias-anisotropic} and proof of Proposition~\ref{prop:bias-isotropic} given in Section~\ref{sec:proofs1}.

\begin{proposition}\label{thm:UB-bias}
Let \(s>0\) and \(L>0\). Let $h\in(0,1/2)^d$ and \(W\in\mathcal{W}^d\) such that for all $i\in\{1,\ldots,d\}$ $W_i$ is of order $\ms$. Then we have:
\begin{equation*}
	\sup_{f\in \tilde{\sobolev}_{s,p}(L)}\|\efn \tilde{f}_{W,h} - f\|_p \leq C_1\|h\|_2^s
\end{equation*}
where $C_1$ is a positive constant that depends only on $W$, $p$ , \(s\) and $L$.

Let $s\in(0,\infty)^d$ and \(L>0\). Let $h\in(0,1/2)^d$ and \(W\in\mathcal{W}^d\) such that for all $i\in\{1,\ldots,d\}$ $W_i$ is of order $\lfloor s_i\rfloor$. Then we have:
\begin{equation*}
	\sup_{f\in {\sobolev}_{s,p}(L)}\|\efn \tilde{f}_{W,h} - f\|_p \leq C_2\sum_{i=1}^dh_i^{s_i}
\end{equation*}
where $C_2$ is a positive constant that depends only on $W$, $p$ , \(s\) and $L$.
\end{proposition}

As we will see in Section~\ref{sec:results}, our boundary kernel estimators and Goldenshluger Lepski selection procedures based on them have also good properties in terms in minimax and adaptive rate of convergence over these classes.


\section{Statistical procedures}\label{sec:stat}

We defined in Section~\ref{ssec:1} a large family of kernels estimators that are well-adapted to the estimation of bounded data. Two subfamilies of estimators designed for the estimation of \emph{isotropic} or \emph{anisotropic} functions are now considered in Sections~\ref{ssec:2} and~\ref{ssec:3} and a unique data-driven procedure is proposed in Section~\ref{ssec:4}.

\subsection{Family of bandwidth and kernels}

We define the set of bandwidth vectors $$\Hn =
\{ h=(h_1,\dotsc,h_d) \in (0,h_n^*]^d : nV_h\geq (\log n)^{c}\}$$ with $c>0$, $h_n^*=\exp(-\sqrt{\log n})$ and $V_h=\prod_{i=1}^d h_i$.

The family of bandwidth $\mathcal{H}_n$ includes in particular for $n$ large enough all the bandwidths $h=(h_1,\ldots,h_d)$ of the form $h_i=n^{-a_i}$ with $0<a_i<1$ and $\sum_{i=1}^d a_i <1$. This family is then rich enough to attain all the optimal rates of convergence of the form $n^{-\overline{s}/(2\overline{s}+1)}$ for $(s_1,\ldots,s_d)\in(0,\infty)^d$ and $\overline{s}=\left(\sum_{i=1}^d\frac{1}{s_i}\right)^{-1}$.
It is possible to have a weaker condition on $h_n^*$ choosing $h_n^*=(\log n)^{-a(p)}$ with $a(p)$ a positive constant that depends on $p$. For the sake of simplicity,  we choose to use $h_n^*=\exp(-\sqrt{\log n})$ to avoid multiple cases in terms of $p$.\\

We consider the family of kernel $(w_m)_{m\in\N}$ defined by:
\begin{equation}\label{eq:omega-m-new}
  w_m(u) = \sum_{r=0}^m \varphi_r(0)\varphi_r(u), \qquad u\in\domain1.
\end{equation}
where $\varphi_k(u)=\sqrt{2k+1}Q_k(2u-1)$ and $Q_k$ is the Legendre Polynomial of degree $k$ on $[-1,1]$ (See \citet{tsybakov2009}). The kernels $w_m$ satisfy several properties given in the following lemma.

\begin{lemma}\label{lem:norme-omega-m}
Set $m\in\N^*$.
The kernel $w_m$ is of order $m$, satisfies $\|{w}_m\|_2=(m+1)$ and
  \begin{equation}\label{eq:2}
    w_m = \argmin_{w\in\cW(m)} \|{w}\|_2
  \end{equation}
  where $\cW(m)\subseteq\cW$ is the family of kernels of order $m$.
  Moreover we have
  \begin{equation}
  \|{w_m}\|_\infty =
    (m+1)^{2}
 \end{equation}
and \begin{equation}\label{eq:omega-m}
  w_m(u) = \sum_{r=0}^m a_r^{(m)} u^r, \qquad u\in\domain1,
\end{equation}
where $a^{(m)}= H_m^{-1} e_0^{(m)}$ where $e_0^{(m)}=(1\, 0 \ldots 0)^\top\in\R^{m+1}$ and $H_{m}=(1/(i+j-1))_{1\le i,j\le m+1}$ is the Hilbert matrix of order $m+1$.\\
\end{lemma}

\begin{figure}[h!]
    \begin{center}
        \includegraphics[width=6cm]{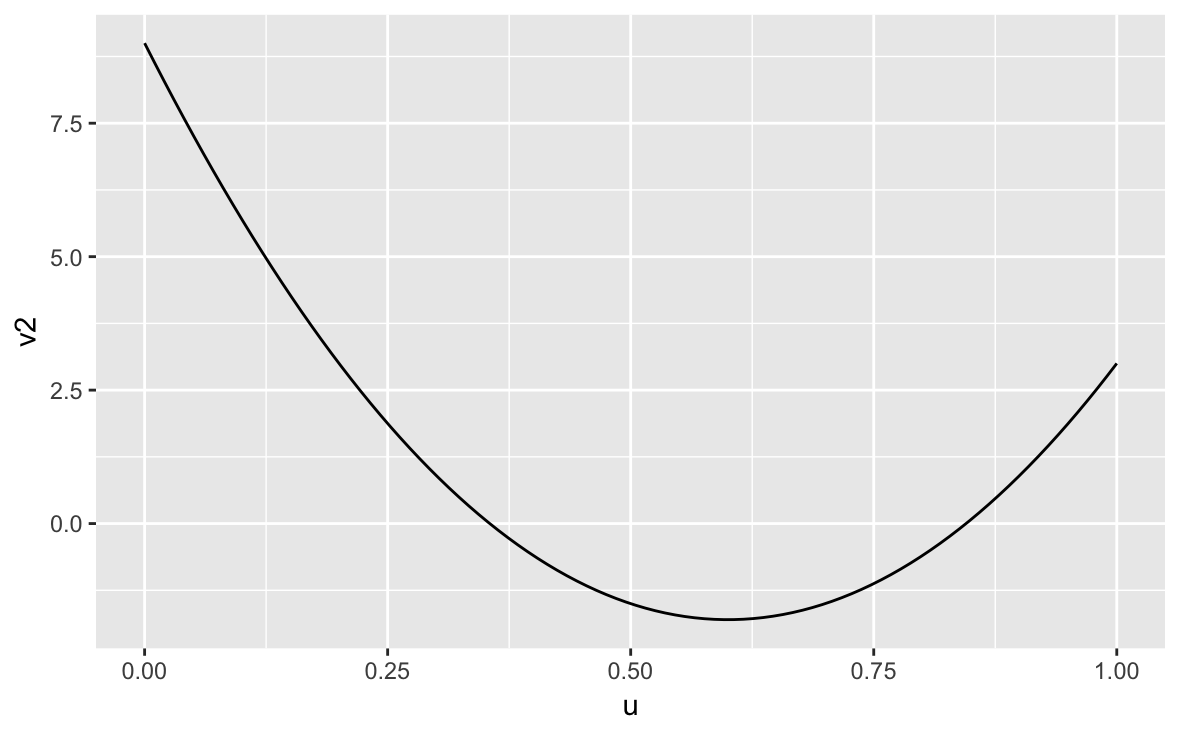}
        \quad
        \includegraphics[width=6cm]{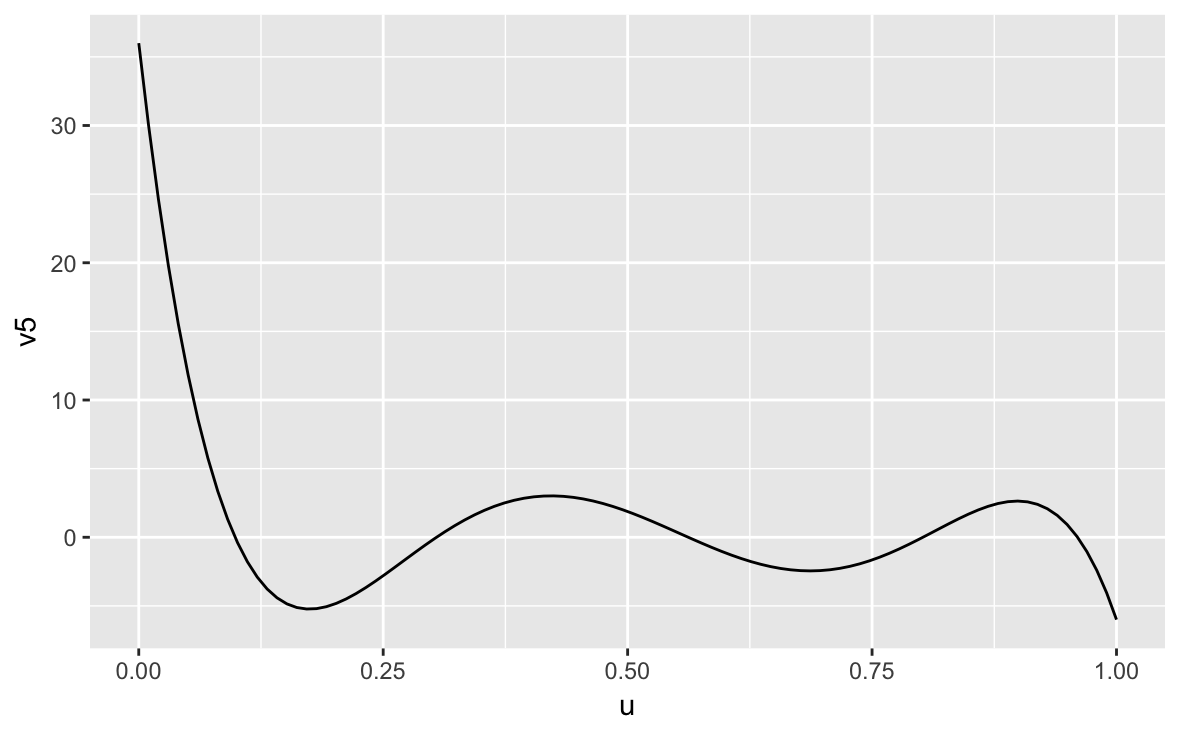}
    \end{center}
    \caption{Plots of the kernel {$w_2$ (left) and $w_5$} (right).}
    \label{fig:kernels}
\end{figure}

Figure~\ref{fig:kernels} represents the kernels $w_m$ for different values of $m$.

\subsection{Isotropic family of estimators}\label{ssec:2}

For $\ell\in\N^*$, we define:
\begin{equation}\label{cond:ml}
  h(\ell)=(e^{-\ell},\ldots,e^{-\ell})
  \quad
  \text{and}
  \quad
  m(\ell) = \left[ \frac{\log n}{2\ell}+\frac12\right],
\end{equation}
where $[ b ]$ stands for the integer part of $b$. We define $\mathcal L_{\mathrm{iso}} = \{ \ell\in\N : h(\ell)\in\Hn\}.$ For any $\ell\in\mathcal L_{\mathrm{iso}}$, we consider $W(\ell)=(w_{m(\ell)},\ldots,w_{m(\ell)})\in\mathcal W^d$ where the univariate kernel $w_m$ is defined by \eqref{eq:omega-m}.

We define the family of estimators $\{\hat f_\ell^{\mathrm{iso}} : \ell\in\mathcal L_{\mathrm{iso}}\}$ where
\begin{equation}
  \hat f_\ell^{\mathrm{iso}} = \tilde f_{W(\ell), h(\ell)}.
\end{equation}

The family $\{\hat f_\ell^{\mathrm{iso}} : \ell\in\mathcal L_{\mathrm{iso}}\}$ contains kernel density estimators constructed with different kernels and bandwidths. Selecting $\ell\in\mathcal{L}_{iso}$, or the estimator $\hat f_\ell^{\mathrm{iso}}$ in this family consists in fact in selecting jointly and automatically the order and the bandwidth of the estimator. The main idea that leads to this construction is the following: if we consider $\ell\approx \log n/(2s+d)$, then $h(\ell)\approx n^{-1/(2s+d)}$ and $m(\ell)\ge s$. In other words, the estimator $\hat f_\ell^{\mathrm{iso}}$ is constructed using a kernel of order greater than $s$ and the usual bandwidth (that is, of the classical order) used to estimate functions with smoothness parameter $s$.
The construction of such a class of estimators allows us to obtain adaptive estimators without any restriction on the smoothness parameter  (see Theorem~\ref{thm:adaptive-isotropic}). However, arbitrary kernels of order $m$ cannot be used to prove Theorem~\ref{thm:oracle2} since a control of the $\mathbb{L}_p$-norm of the kernels is required. In particular in Lemma~\ref{lem:norme-omega-m}, we give bounds  on the $\mathbb{L}_p$-norm of $w_m$ and we prove that $w_m$ is the kernel of order $m$ with the smallest $L_2$ norm within the kernels of $\mathcal{W}$ of order $m$.

Note that simultaneous choice of kernel and bandwidth has already been used in the framework of sharp adaptive estimation only for pointwise and sup-norm risks. On the one hand, in the Gaussian white noise model, \citet{bertin2005} and \citet{lepskispokoiny97} assume that the smoothness parameter is less than or equal to $2$. On the other hand, in the density model, \citet{butucea2001} consider the case of a finite grid of integer smoothness parameters $s_1<\ldots<s_R$ and propose an adaptive procedure for the pointwise risk over a scale of classical Sobolev classes. Note that in this paper, the maximal smoothness parameter $s_R$ may tend to infinity as $n$ goes to infinity. To our understanding this possibility relies on the fact that the kernels are uniformly bounded by a constant that depends only on $s_1$. Studying the risk over classical Sobolev classes on $\R$ allows \citet{butucea2001} to define the kernels on the Fourier domain and to replace the moment condition \eqref{eq:moment-condition} by a weaker one \citep[see][Section~1.3 for more details]{tsybakov2009}.

Our framework is very different. Indeed we consider the estimation in $\mathbb{L}_p$ risks of densities with compact support that belongs to a scale of Sobolev-Slobodetskii classes indexed by a smoothness parameter $s\in\R_+$. To do so we have to consider the classical moment condition~\eqref{eq:moment-condition} which implies, according to Lemma~\ref{lem:norme-omega-m}, that the sup-norm of any kernel of order $m$ tends to infinity with $m$. This requires more technical control of the stochastic terms to obtain the minimax rates of convergence without additional logarithmic factor.

\subsection{Anisotropic family of estimators}\label{ssec:3}
\label{sec:anisotropic}
Let $W^\circ=(W^\circ_1,\ldots, W^\circ_d)\in\mathcal W^d$ be such that, for any $i=1,\dotsc,d$, $W_i^\circ$ is a bounded kernel and consider $h(\ell)=(h_1(\ell),\ldots, h_d(\ell))$ defined by:
\begin{equation}\label{eq:def-hl}
  h_i(\ell) = e^{-\ell_i},
  \qquad
  i=1,\ldots,d
\end{equation}
where $\ell\in\mathcal L_{\mathrm{ani}}=\{\ell\in\N^d : h(\ell)\in\Hn\}$.

We define the anisotropic family of estimators $\{\hat f_\ell^{\mathrm{ani}} : \ell\in\mathcal L_{\mathrm{ani}}\}$ where
\begin{equation}
  \hat f_\ell^{\mathrm{ani}} = \tilde f_{W^\circ, h(\ell)}.
\end{equation}

To make the notation similar to the isotropic case we define $W(\ell)=W^\circ, \forall\ell\in\mathcal L_{\mathrm{ani}}$. Note that this family of estimators is more classical than the one constructed in the previous section. All the estimators are defined using the same kernel $W^\circ$ and depend only on a multivariate bandwidth. Nevertheless, in the following (see Theorem~\ref{thm:adaptive-anisotropic}), we will choose a kernel $W^{\circ}=(W^\circ_1,\ldots,W^\circ_d)$  such that for all $i\in\{1,\ldots,d\}$ $W^\circ_i$ is of order $M_i$ and a possible candidate is $W^\circ_i=w_{M_i}$.

\subsection{Selection rule}\label{ssec:4}
\label{sec:selection-rule}
Although the two families differ, the selection procedure is the same in both cases. For the sake of generality, we introduce the following notation: $\mathcal L$ is either $\mathcal L_{\mathrm{ani}}$ or $\mathcal L_{\mathrm{iso}}$ and $\hat f_\ell$ then denotes $\hat f_\ell^{\mathrm{ani}}$ or $\hat f_\ell^{\mathrm{iso}}$. For $\eps\in\{0,1\}^d$, $h\in\Hn$ and $W\in\mathcal W^d$ we define:
\begin{equation}
  \Delta_{d,\varepsilon} = \prod_{i=1}^d \left(\frac{\varepsilon_i}2, \frac{1+\varepsilon_i}2\right),\qquad \|W\|_p = \left\|\bigotimes_{i=1}^d W_i\right\|_p
\end{equation}
\begin{align}\label{eq:hatlambda}
  \hat\Lambda_\varepsilon(W,h,p)
  &= \sqrt{V_{h}}\left(\int_{\domain{d,\varepsilon}}\left(\frac1n\sum_{j=1}^n \mathcal K_{W,h}^2(t,X_j)\right)^{p/2} \D t\right)^{1/p}
\end{align}
and
\begin{equation}
  \hat\Gamma_\varepsilon(W,h,p) = \left\{
  \begin{array}{cc}
    2^{-\frac{d(2-p)}{2p}} \|W\|_2 & \text{if $1\leq p\leq2$}\\
    C_p^* \left(\hat\Lambda_\varepsilon(W,h,p) + 2\|W\|_p\right) & \text{if $p>2$}\label{eq:hatgamma}
  \end{array}
  \right.
\end{equation}
where $C_p^* = {14.7 p}/{\log p}$ is the best known constant in the Rosenthal inequality \citep[see][]{bestconstant}. For any $\ell,\ellp\in\mathcal L$ we consider:

\begin{equation}
  \widehat{M}_p(\ell) = \frac{1}{\sqrt{nV_{h(\ell)}}} \sum_{\varepsilon\in\{0,1\}^d} \hat\Gamma_\varepsilon(W(\ell), h(\ell), p)
\quad
\text{and}
\quad
  \hat M_p(\ell,\ellp) =
  \hat M_p(\ellp) + \hat M_p(\ellp\wedge\ell)
\end{equation}
where $\ell\wedge\ellp$ is the vector with coordinates $\ell_i\wedge\ell_i'=\min(\ell_i,\ell_i')$. Now, for any $\tau>0$ we define:
\begin{align}\label{eq:eqb}
  \hat B_p(\ell)
  &= \max_{\ellp\in\mathcal L} \left\{\|{\hat f_{\ell\wedge\ellp}-\hat f_{\ellp}}\|_p-(1+\tau)\hat M_p(\ell,\ellp)\right\}_+
\end{align}
where $x_{+}=\max(x,0)$ denotes the positive part of $x$.

We then select
\begin{equation}
  \hat\ell = \argmin_{\ell\in\mathcal L} \left(\hat B_p(\ell) + (1+\tau)\hat M_p(\ell)\right)
\end{equation}
which leads to the final plug-in estimator defined by $\hat f = \hat f_{\hat\ell}$.
In what follows we denote by $\hat f^{\mathrm{ani}}$ and $\hat f^{\mathrm{iso}}$ the resulting estimators.

\begin{remark}
    {This procedure is inspired by the method developed by Goldenshluger and Lepski. Here $\hat B_p(\ell)$ is linked with the bias term of the estimator $\hat f_\ell$, see~\eqref{eq:par1}, and $\hat M_p(\ell)$ is an empirical version of an upper bound on the standard deviation of this estimator.
    In fact, for $p\le 2$, the standard deviation in $\mathbb{L}_p$-norm of $\hat f_\ell$ on $\Delta_{d,\varepsilon}$ is bounded by $2^{-\frac{d(2-p)}{2p}}\|W(\ell)\|_2$. For $p>2$, the bound depends on $f$ (see Lemma~\ref{lem:espY}), that is the reason why we use an empirical version of this bound defined in \eqref{eq:hatgamma}.
    This implies that $\hat f$ realizes a trade-off between $\hat B_p(\ell)$ and $(1+\tau)\hat M_p(\ell)$. This can be interpreted as an empirical counterpart of the classical trade-off between the bias and the standard deviation. Note that as discussed in \citet{lacour2016minimal} it is also possible to consider in \eqref{eq:eqb} a different constant $\tau'$ satisfying $\tau'<\tau$.}

\end{remark}

\section{Results}\label{sec:results}


In this section we present our results. Theorem~\ref{thm:oracle} consists in an oracle-type inequality which guarantees that the anisotropic estimation procedure defined above performs almost as well as the best estimator from the collection $\{\hat f_\ell^\mathrm{ani} : \ell\in\mathcal L_{\mathrm{ani}}\}$. Moreover, Theorem~\ref{thm:adaptive-anisotropic} states that this procedure also achieves the minimax rate of convergence simultaneously over each anisotropic Sobolev-Slobodetskii class in a given scale.

\begin{theorem}\label{thm:oracle}
  Assume that  $f:\domain{d}\to\R$ is a density function such that $\|f\|_\infty\leq F_\infty$. Then there exists a positive constant $\const1$ that depends only on $F_\infty$, $W^\circ$, $p$, $q$ and $\tau$, such that, for any $n\geq2$:
  \begin{equation}
    R_n^{(p,q)}(\hat f^{\mathrm{ani}},f) \leq \const1 \inf_{\ell\in\mathcal L_{\mathrm{ani}}} \left\{  \| \e_f^n \hat{f}_{\ell}^{\mathrm{ani}}-f \|_p
    + \max_{\ellp\in\mathcal L_{\mathrm{ani}}} \| \e_f^n \hat{f}_{\ellp}^{\mathrm{ani}}-\e_f^n \hat{f}_{\ell\wedge\ellp}^{\mathrm{ani}} \|_p
    + \frac{1}{(nV_{\hell})^{1/2}}\right\}.
  \end{equation}
\end{theorem}
\begin{theorem}\label{thm:adaptive-anisotropic}
  Set $M=(M_1,\ldots,M_d)\in\N^d$, $s\in\prod_{i=1}^d(0,M_i+1]$ and $L>0$. Assume that $W^\circ$ is such that $W^\circ_i$ is of order greater than or equal to $M_i$. Then, our estimation procedure $\hat f^{\mathrm{ani}}$ is such that:
  \begin{equation}
    \limsup_{n\to+\infty} n^{\frac{\bar s}{2\bar s+1}}\sup_{f \in\sobolev_{s,p}(L)} R_n^{(p,q)}(\hat f^{\mathrm{ani}},f) < +\infty.
  \end{equation}
  Moreover, if $s=(s_1,\dotsc,s_d)$ is such that any $s_i$ is not an integer, the following property is satisfied:
  \[
  	\liminf_{n\to\infty}  n^{\frac{\bar s}{2\bar s+1}} \inf_{\tilde f} \sup_{f\in \sobolev_{s,p}(L)} R_n^{(p,q)}(\tilde f,f) >0
  \]
  where the infimum is taken over all possible estimators.
\end{theorem}
Theorems~\ref{thm:oracle2} and~\ref{thm:adaptive-isotropic} are the analogues of Theorems~\ref{thm:oracle} and~\ref{thm:adaptive-anisotropic} respectively, transposed to the isotropic estimation procedure. Note however that the scale of functional classes considered in Theorem~\ref{thm:adaptive-isotropic} is huge since there is no restriction on the smoothness parameter $s>0$, contrary to classical results (including Theorem~\ref{thm:adaptive-anisotropic}).

\begin{theorem}\label{thm:oracle2}
  Assume that  $f:\domain{d}\to\R$ is a density function such that $\|f\|_\infty\leq F_\infty$. Then there exists a positive constant $\const2$ that depends only on $F_\infty$, $p$, $q$ and $\tau$, such that, for any $n\geq2$:
  \begin{equation}
    R_n^{(p,q)}(\hat f^{\mathrm{iso}},f) \leq \const2 \inf_{\ell\in\mathcal L_{\mathrm{iso}}} \left\{ \max_{\ellp\geq \ell}  \| \e_f^n \hat{f}_{\ellp}^{\mathrm{iso}}-f \|_p    +  \frac{\|W(\ell)\|_{p\vee 2}}{(nV_{\hell})^{1/2}} \right\}.
  \end{equation}
\end{theorem}
\begin{theorem}\label{thm:adaptive-isotropic}
    Set $s>0$ and $L>0$. We have:
    \begin{equation}
        \limsup_{n\to+\infty} n^{\frac{s}{2s+d}}\sup_{f \in\tilde\sobolev_{s,p}(L)} R_n^{(p,q)}(\hat f^{\mathrm{iso}},f) < +\infty
    \end{equation}
    and, if $s$ is not an integer,
    \begin{equation}
    	\liminf_{n\to\infty}  n^{\frac{s}{2s+1}} \inf_{\tilde f} \sup_{f\in\tilde\sobolev_{s,p}(L)} R_n^{(p,q)}(\tilde f,f) >0
    \end{equation}
    where the infimum is taken again over all possible estimators.

\end{theorem}
%
\begin{remark}\label{rem:main-results}
	Theorems~\ref{thm:adaptive-anisotropic} and~\ref{thm:adaptive-isotropic} are established for scales of Sobolev-Slobodetskii classes. However similar results are still true if one replaces these classes with classical (an)isotropic Hölder classes. Remark also that the lower bounds are proved for non-integer smoothness parameters. As mentioned above, in this situation, the Sobolev-Slobodetskii classes correspond to usual Besov spaces.

  In Theorems~\ref{thm:oracle} and~\ref{thm:oracle2}, the right hand sides of the equations can be easily interpreted.  In both situations, the term $(nV_{\hell})^{-1/2}$ is of the order of the standard deviation of $\hat f_\ell$. Moreover the terms $\max_{\ellp\in\mathcal L_{\mathrm{ani}}} \| \e_f^n \hat{f}_{\ellp}^{\mathrm{ani}}-\e_f^n \hat{f}_{\ell\wedge\ellp}^{\mathrm{ani}} \|_p$ and $\max_{\ellp\geq \ell}  \| \e_f^n \hat{f}_{\ellp}^{\mathrm{iso}}-f \|_p$ are  linked with the bias of this estimator. More precisely, Proposition~\ref{prop:bias-anisotropic} and Proposition~\ref{prop:bias-isotropic} ensure that these terms  have the same behaviour as the bias term
  $\| \e_f^n \hat{f}_{\ell}-f \|_p$ as soon as $f$ belongs to Sobolev-Slobodetskii classes.

  Finally, to our best knowledge, even in the case of density with support in $\R$, adaptive results in $\mathbb{L}_p$ without restriction on the smoothness parameter  as in Theorem~\ref{thm:adaptive-isotropic} are not known for either the Sobolev-Slobodetskii classes or the H\"older classes. This is not the case in Theorem~\ref{thm:adaptive-anisotropic} where the adaptive result is obtained only for $s\in\prod_{i=1}^d(0,M_i+1]$ where the $M_i$ are the orders of the kernel $W^\circ$.  The main difference between the isotropic case and the anisotropic case lies in the control of the quantity  $B_p(\ell)$ which is linked with the terms $\|\mathbb{E}\hat f_{\ell'}- \mathbb{E}\hat f_{\ell\wedge\ell'}\|_p$ for $\ell'\in\mathcal{L}$. In the isotropic case, if $\ell'\leq\ell$, these terms vanish and it remains to control
  \begin{equation}\label{XX}\max_{\ell'\geq\ell}\|\mathbb{E}\hat f_{\ell'}- \mathbb{E}\hat f_{\ell\wedge\ell'}\|_p\le 2\max_{\ell'\geq\ell}\|\mathbb{E}\hat f_{\ell'}- f\|_p.
  \end{equation}
  The study of \eqref{XX} involves Taylor expansion of $f$ and each estimator  $\hat{f}_\ell$ can be based on a different kernel. In the anisotropic case, \eqref{XX} is never more valid  and $\|\mathbb{E}\hat f_{\ell'}- \mathbb{E}\hat f_{\ell\wedge\ell'}\|_p$  can be expressed in terms of  the difference of $f$ in two different values (in order to use a Taylor expansion) only when $\hat{f}_{\ell'}$ and $\hat f_{\ell\wedge\ell'}$ are based on the same kernel.
\end{remark}

\section{Proofs}\label{sec:proofs1}

The proofs of Theorems~\ref{thm:oracle}--\ref{thm:adaptive-isotropic} are based on propositions and lemmas which are given below. Before stating these results, we introduce some notation that are used throughout the rest of the paper. For $W=(W_1,\ldots,W_d)\in \mathcal{W}^d$, $h\in\Hn$ and $\varepsilon\in\{0,1\}^d$, we define the quantity:
\begin{equation}
  \Gamma_\varepsilon(W,h,p) = \left\{\begin{array}{cc}
 2^{-\frac{d(2-p)}{2p}} \|W\|_2 & \text{if $1\leq p\leq2$}\\
 C_p^* \left(\Lambda_\varepsilon(W,h,p) + 2\|W\|_p\right) & \text{if $p>2$}
  \end{array}\right.
\end{equation}
where
\begin{align}
   \Lambda_\varepsilon(W,h,p)
   &= \sqrt{V_h}\left(\int_{\domain{d,\varepsilon}} \left(\int_{\domain{d}} \mathcal K_{W,h}^2(t,x) f(x) \D x\right)^{p/2} \D t\right)^{1/p}.
\end{align}

For $g:\domain{d}\to\R$ and $r\ge 1$ we denote
\begin{equation*}
\|g\|_{r,\varepsilon}=\left(\int_{\Delta_{d,\varepsilon}}|g(x)|^r \D x\right)^{1/r}.
\end{equation*}

The process $\xi_{W,h}$ is defined by
\begin{equation*}
\xi_{W,h}(t)= \left(\frac{V_h}{n}\right)^{1/2}\sum_{j=1}^n \left( \mathcal{K}_{W,h}(t,X_j)-\e_f^n \mathcal{K}_{W,h}(t,X_j)\right), \qquad t\in\Delta_d.
\end{equation*}

Finally, for $\ell\in\mathcal L$ we define (using the generic notation for the isotropic and the anisotropic cases):
\begin{equation}
W^*(\ell)=\left(\frac{(W_1(\ell))^2}{\|W_1(\ell)\|_2^2},\ldots,\frac{(W_d(\ell))^2}{\|W_d(\ell)\|_2^2}\right).
\end{equation}


\begin{proposition}[Anisotropic case]\label{prop:bias-anisotropic}
    Set $M=(M_1,\ldots,M_d)\in\N^d$. Assume that $W^\circ$ is such that $W^\circ_i$ is of order greater than or equal to $M_i$. Set $s=(s_1,\ldots, s_d)\in\prod_{i=1}^d(0,M_i]$  and $L>0$. Then, for any $f\in\sobolev_{s,p}(L)$:
    \begin{equation}\label{eq:bias1-anisotropic}
        \|\efn \hat f_\ell^{\mathrm{ani}} - f\|_p \leq 2^{d/p} d \biggl(\prod_{i=1}^d (M_i+1)\biggr) L \sum_{i=1}^d \big(h_i(\ell)\big)^{s_i}.
    \end{equation}
    \begin{equation}\label{eq:bias2-anisotropic}
        \max_{k\in\mathcal L_{\mathrm{ani}}}\|\efn \hat f_{k}^{\mathrm{ani}} - \efn \hat f_{\ell\wedge k}^{\mathrm{ani}}\|_p \leq 2^{1+d/p} d\left(\prod_{i=1}^d (M_i+1)\right)L\sum_{i=1}^d \big(h_i(\ell)\big)^{s_i}.
    \end{equation}
\end{proposition}

\begin{proposition}[Isotropic case]\label{prop:bias-isotropic}
    Set $s>0$  and $L>0$. Then for any $\ell\in\mathcal L_{\mathrm{iso}}$ we have:
    \begin{equation}\label{eq:bias1-isotropic}
        \sup_{f\in\tilde\sobolev_{s,p}(L)}\max_{\ellp\geq\ell}\|\efn \hat f_{\ellp}^{\mathrm{iso}} - f\|_p
        \leq
        \const{3}\left(\|W(\ell)\|_{\infty}L(h_1(\ell))^s
        + \sqrt{\frac{h_n^*}{n}}\right),
    \end{equation}
    where the positive constant $\const{3}$ depends only on $d$, $p$, $s$ and $L$.
\end{proposition}

\begin{proposition}\label{prop:majorant}
Set $p,q\geq1$. Assume that $f$ is such that $\|f\|_\infty\le F_\infty$.
\begin{itemize}
\item Let $\ell\in\mathcal{L}_{\mathrm{iso}}$. There exists a positive constant $\const4$ that depends only on $p$, $q$, $\tau$  and $F_\infty$ such that
\begin{equation}
    \e_f^n\left\{\|\hat f_{\ell}^{\mathrm{iso}}-\e_f^n \hat f_{\ell}^{\mathrm{iso}}\|_p-(1+\tau)\widehat{M}_p(\ell)\right\}_+^q
\leq\const4 n^{-q}.
\end{equation}
\item Let $\ell\in\mathcal{L}_{\mathrm{ani}}$. There exists a positive constant $\const5$ that depends only on $p$, $q$, $\tau$, $W^\circ$ and $F_\infty$ such that
\begin{equation}
    \e_f^n\left\{\|\hat f_{\ell}^{\mathrm{ani}}-\e_f^n \hat f_{\ell}^{\mathrm{ani}}\|_p-(1+\tau)\widehat{M}_p(\ell)\right\}_+^q
\leq\const5 n^{-q}.
\end{equation}
\end{itemize}
\end{proposition}

\begin{lemma}\label{lem:espY}
 Assume that $f$ satisfies $\|f\|_\infty\leq \mathrm F_\infty$. For any $W\in\mathcal{W}^d$, $r\geq 1$ and $h\in \Hn$, we have:
     \begin{equation}
         \e_f^n \|\xi_{W,h}\|_{r,\varepsilon} \leq \Gamma_\varepsilon(W,h,r)\le C_0\|W\|_{2\vee r},
     \end{equation}
     where $C_0$ is an absolute constant that depends only on $r$ and $\mathrm F_\infty$.
 \end{lemma}

\begin{lemma}\label{lem:bousq}Assume that $f$ satisfies $\|f\|_\infty\leq \mathrm F_\infty$.
     For any $W\in\mathcal{W}^d$, $r\geq 1$ and $h\in \Hn$, we have:
     \begin{align}\label{eq:11}
     \p(\|\xi_{W,h}\|_{r,\varepsilon}\hspace{-0.1cm}-\hspace{-0.1cm}\e_f^n \|\xi_{W,h}\|_{r,\varepsilon} \geq \frac{\tau}{2}\Gamma_\varepsilon(W,h,r)\hspace{-0.1cm}+\hspace{-0.1cm}x)
     &\leq \exp\left(-\frac{C_2x^2(\alpha_n(r))^{-1}}{\|W\|^2_{2\vee r}+x\|W\|_{r}}\right) \exp \left(-C_1 \alpha_n(r)\right)
 \end{align}
 where $C_1$ and $C_2$ are absolute constants that depend only on $r$, $\tau$ and $\mathrm F_\infty$,
 \begin{equation}
   \alpha_n(r) =
   \left\{\begin{array}{cc}
     (h_n^*)^{-d\left(\frac2r-1\right)} & \text{if $1\leq r<2$} \\
     (h_n^*)^{-\frac dr} & \text{if $r\geq 2$}
   \end{array}\right.
 \end{equation}
 and $h_n^*=\exp(\sqrt{-\log n})$.
 \end{lemma}

We finally state the following lemma that allows us to bound the bias terms which appear in the oracle inequality.

\begin{lemma}\label{lem:bias1}
    Let $h=(h_1,\ldots,h_d)$ and $\eta=(\eta_1,\ldots,\eta_d)$ be two bandwidths in $\mathcal H_{n}$ such that $\eta_i\in\{0,h_i\}$. Set $W=(w_{M_1},\ldots,w_{M_d})\in\mathcal W^d$ and define:
    \begin{equation}
    S_{W,h, \eta}^*(f) =
        \left(
        \int_{\Delta_{d,\mathbf 0}}
        \left\vert \int_{\Delta_d}
        \left(\prod_{i=1}^d w_{M_i}(u_i)\right)
        [f(t+h\cdot u)-f(t+\eta\cdot u)] \D u \right\vert^p \D t
        \right)^{1/p}
    \end{equation}
    where $h\cdot u$ denotes the coordinate-wise product of the vectors $h$ and $u$. Assume that $f$ belongs to $\sobolev_{s,p}(L)$ and that, for any $i=1,\ldots,d$, the kernel $W_i$ is of order greater than or equal to $\lfloor s_i\rfloor$.
Then we have:
    \begin{equation}
        S_{W,h, \eta}^*(f) \leq d\left(\prod_{i=1}^d (M_i+1)\right)L \sum_{i\in I} h_i^{s_i}
    \end{equation}
    where $I=\{i=1,\ldots,d : \eta_i=0\}$.
\end{lemma}

\subsection{Proof of Proposition~\ref{thm:LB-bias}}

We note that:
\begin{equation}\label{eq:thm1-1}
	\int_{\Delta_1} \big|\e_{f_0}^n\hat{f}_h(t)-f_0(t)\big|^p dt
	\geq
	\int_{0}^{\delta h} \big|\e_{f_0}^n\hat{f}_h(t)-1\big|^p dt.
\end{equation}
Now we remark that, for any $t\in(0,\delta h)$, we have $(t-1)/h\leq -1$ which implies that:
\begin{align}
	\e_{f_0}^n\hat{f}_h(t)
	&= \int_\R K_h(t-u)\1_{(0,1)}(u) du \\
	&= \int_{(t-1)/h}^{t/h} K(u) du \\
	&= \int_{-1}^{t/h} K(u) du.
\end{align}
Since $t/h\leq \delta$ we obtain that in this situation $\e_{f_0}^n\hat{f}_h(t)\leq 1-\gamma$. As a consequence, for any $t\in(0,\delta h)$ we have:
\[
	f_0(t)-\e_{f_0}^n\hat{f}_h(t) \geq \gamma.
\]
Combining this inequality with~\eqref{eq:thm1-1} we obtain:
\[
	\int_{\Delta_1} \big|\e_{f_0}^n\hat{f}_h(t)-f_0(t)\big|^p dt
	\geq
	\int_0^{\delta h} \gamma^p dt.
\]
Proposition~\ref{thm:LB-bias} follows.

\subsection{Proof of Proposition~\ref{thm:LB-risk}}

Let $f\in\Sigma$ be a density function and let $h\in(0,1/4)$. Using Jensen inequality we obtain for any $t\in\Delta_1$:
\begin{equation}\label{eq:thm2-1}
	|\e_f^n\hat{f}_h(t)-f(t) |^p
	\leq
	\e_f^n \left|\hat{f}_h(t)-f(t)\right|^p.
\end{equation}
Integrating over $\Delta_1$ we obtain:
\begin{equation}\label{eq:thm2-2}
	\|\e_f^n\hat{f}_h-f \|_p
	\leq
	R_n^{(p,p)}(\hat{f}_h,f).
\end{equation}
Now, using the triangular inequality we have:
\begin{equation}
	\|\hat{f}_h-\e_f^n \hat{f}_h\|_p
	\leq
	\|\hat{f}_h-f\|_p+\|\e_f^n \hat{f}_h-f\|_p.
\end{equation}
Using again the triangular inequality we obtain:
\begin{equation}\label{eq:thm2-3}
	\left(\e_f^n\|\hat{f}_h-\e_f^n \hat{f}_h\|_p^p\right)^{1/p}
	\leq R_n^{(p,p)}(\hat{f}_h,f) + \|\e_f^n\hat{f}_h-f \|_p.
\end{equation}
Combining~\eqref{eq:thm2-2} and~\eqref{eq:thm2-3} we obtain:
\begin{equation}
	3 R_n^{(p,p)}(\hat{f}_h,f) \geq \|\e_f^n\hat{f}_h-f \|_p + \left(\e_f^n\|\hat{f}_h-\e_f^n \hat{f}_h\|_p^p\right)^{1/p}.
\end{equation}
Fixing $f=f_0$ and using Theorem~\ref{thm:LB-bias} (note that $h\leq1/4$ implies that $h\leq(1+\delta)^{-1}$) we obtain:
\begin{equation}\label{eq:thm2-10}
	3 R_n^{(p,p)}(\hat{f}_h,f_0)
	\geq
	(\gamma \delta^{1/p}) h^{1/p} + \left(\e_{f_0}^n\|\hat{f}_h-\e_{f_0}^n \hat{f}_h\|_p^p\right)^{1/p}.
\end{equation}
Now, it remains to bound the last term of the right hand side of this inequality. To this aim note that:
\begin{align}
	\e_{f_0}^n\|\hat{f}_h-\e_{f_0}^n \hat{f}_h\|_p^p
	&\geq
	\e_{f_0}^n \int_h^{1-h}|\hat{f}_h(t)-\e_{f_0}^n \hat{f}_h(t)|^p dt\\
	&\geq
	\int_h^{1-h}\left(\e_{f_0}^n |\hat{f}_h(t)-\e_{f_0}^n \hat{f}_h(t)|^2\right)^{p/2} dt
\end{align}
where the last line follows from Jensen inequality. Using that $t\in(h,1-h)$ and that $\operatorname{Supp}(K)\subseteq[-1,1]$ we obtain:
\[
	\e_{f_0}^n \hat{f}_h(t) = 1.
\]
This implies that for any $t\in(h,1-h)$:
\begin{align}
	\e_{f_0}^n |\hat{f}_h(t)-\e_{f_0}^n \hat{f}_h(t)|^2
	&= \frac1n \left(\e_{f_0}K_h^2(t-X)-1\right)\\
	&= \frac1{nh} \int_\R K^2(u)\1_{(0,1)}(t-hu) du -\frac1n\\
	&= \frac1{nh} \int_{-1}^1 K^2(u) du -\frac1n\\
	&\geq \frac{\|K\|_2^2}{2nh}.
\end{align}
Last inequality holds since $h\leq 1/4\leq\|K\|_2^2/2$ (using Cauchy-Schwarz inequality). Finally we obtain:
\begin{align}
	\left(\e_{f_0}^n\|\hat{f}_h-\e_{f_0}^n \hat{f}_h\|_p^p\right)^{1/p}
	&\geq \left(\frac{(1-2h)\|K\|_2^2}{2nh}\right)^{1/2}\\
	&\geq \frac{\|K\|_2}{2} \, (nh)^{-1/2}\label{eq:thm2-11}
\end{align}
Last inequality holds since $h\leq1/4$. Now, combining~\eqref{eq:thm2-10} with~\eqref{eq:thm2-11} and minimizing with respect to $h$, Proposition~\ref{thm:LB-risk} follows.

\subsection{Proof of
Proposition~\ref{prop:bias-anisotropic}}

We first prove~\eqref{eq:bias1-anisotropic}. {Set $W=W^\circ$ and $h=h(\ell)$. We have $\hat f^{\mathrm{ani}}_\ell=\tilde f_{W,h}$ and
\begin{align}
    \|\e_f^n \tilde f_{W,h}-f\|_p^p
    &= \sum_{\varepsilon\in\{0,1\}^d} \int_{\Delta_{d,\varepsilon}} |\e_f^n \tilde f_{W,h}(t)-f(t)|^p \D t\\
    &= \sum_{\varepsilon\in\{0,1\}^d} \int_{\Delta_{d,\varepsilon}}
    \left\vert \int_{\Delta_d} \mathcal K_{W,h}(t,x) f(x) \D x-f(t) \right\vert^p \D t\\
    &=\sum_{\varepsilon\in\{0,1\}^d} \int_{\Delta_{d,\mathbf 0}}
    \left\vert \int_{\Delta_d} \mathcal K_{W,h}(u,y) f_\varepsilon(y) \D y-f_\varepsilon(u) \right\vert^p \D u,\label{eq:symomega}
\end{align}
 where
 \begin{equation}
    f_{\varepsilon}(u) = f(\ldots,u_i(1-\varepsilon_i)+(1-u_i)\varepsilon_i,\ldots).
\end{equation}
Line \eqref{eq:symomega} is obtained doing, for each $\varepsilon\in\{0,1\}^d$, the changes of variables in both integrals, $t_i=u_i(1-\varepsilon_i)+(1-u_i)\varepsilon_i$ and $x_i=y_i(1-\varepsilon_i)+(1-y_i)\varepsilon_i$ for all $i\in\{1,\ldots,d\}$, and using that $K_{W,h}(t,x)=K_{W,h}(u,y)$.
As a consequence
\begin{equation*}
\|\e_f^n\tilde f_{W,h}-f\|_p^p=\sum_{\varepsilon\in\{0,1\}^d} \left(S_{W,h}(f_\varepsilon)\right)^p
\end{equation*}
where
\begin{equation}\label{eq:swh}
S_{W,h}(f) = \left(\int_{\Delta_{d,\mathbf 0}}
    \left\vert \int_{\Delta_d} \mathcal K_{W,h}(t,x) f(x) \D x-f(t) \right\vert^p \D t\right)^{1/p}.
\end{equation}
Since $f\in\sobolev_{s,p}(L)\iff f_{\varepsilon}\in\sobolev_{s,p}(L)$, we obtain
\begin{equation}\label{eq:bias-general}
    \sup_{f\in\sobolev_{s,p}(L)} \|\e_f^n \tilde f_{W,h}-f\|_p\leq  2^{d/p} \sup_{f\in\sobolev_{s,p}(L)} S_{W,h}(f).
\end{equation}
Then Equation~\eqref{eq:bias1-anisotropic} follows from Lemma~\ref{lem:bias1} and the fact that $S_{W,h}(f)=S_{W,h,0}^*(f)$.
}

Now, let us prove~\eqref{eq:bias2-anisotropic}. Set $h=h(k)$ and $h'=h(k\wedge\ell)=h(k)\vee h(\ell)$. Similarly to \eqref{eq:symomega} we have:{
\begin{align}
\sup_{f\in\sobolev_{s,p}(L)} &\|\e_f^n \tilde f_{W,h}-\efn \tilde f_{W,h'}\|_p^p\\
    &\leq 2^d \sup_{f\in\sobolev_{s,p}(L)}\int_{\Delta_{d,\mathbf 0}}
    \left\vert
    \int_{\Delta_d} \mathcal K_{W,h}(t,x) f(x) \D x
    -
    \int_{\Delta_d} \mathcal K_{W,h'}(t,x) f(x) \D x
    \right\vert^p \D t\\
    &\leq 2^d  \sup_{f\in\sobolev_{s,p}(L)}
    \int_{\Delta_{d,\mathbf 0}}
        \left\vert
        \int_{\Delta_d} \left(\prod_{i=1}^d W_i(u_i) \right)
        [f(t+h\cdot u) - f(t+h'\cdot u)]\D u
        \right\vert^p \D t.
\end{align}}
Let $\eta=(\eta_1,\ldots,\eta_d)$ be a bandwidth defined by
\begin{equation}
    \eta_i = \begin{cases}
        0 &\text{if } h_i<h_i'\\
        h_i &\text{if } h_i=h_i'.
    \end{cases}
\end{equation}
We have:{
\begin{align}
     \sup_{f\in\sobolev_{s,p}(L)}&\|\e_f^n \tilde f_{W,h}-\efn \tilde f_{W,h'}\|_p^p\\
  &  \leq
    2^{d+p}  \sup_{f\in\sobolev_{s,p}(L)}\max_{H\in\{h,h'\}}
    \int_{\Delta_{d,\mathbf 0}}
    \left\vert
    \int_{\Delta_d} \left(\prod_{i=1}^d W_i(u_i) \right)
    [f(t+H\cdot u) - f(t+\eta\cdot u)]\D u
    \right\vert^p \D t.
\end{align}}
Using Lemma~\ref{lem:bias1}, we obtain:
\begin{align}
    \sup_{f\in\sobolev_{s,p}(L)} \|\e_f^n \tilde f_{W,h}-\efn \tilde f_{W,h'}\|_p
    \leq
    2^{1+d/p} d\left(\prod_{i=1}^d (M_i+1)\right)L
    \max_{H\in\{h,h'\}} \sum_{i\in I} H_i^{s_i}
\end{align}
where $I=\{i : \eta_i=0\}$. Since $H_i\leq h_i(\ell)$ for any $i\in I$, this allows us to conclude.

\subsection{Proof of Proposition~\ref{prop:bias-isotropic}}

In the same way that the proof of Proposition~\ref{prop:bias-anisotropic}, we obtain:
\begin{equation}
    \sup_{f\in\tilde S_{s,p}(L)} \|\efn \hat f^{\mathrm{iso}}_\ell - f\|_p
    \leq
    2^{d/p} \sup_{f\in\tilde S_{s,p}(L)} S_{W(\ell),h(\ell)}(f),
\end{equation}
where $S_{W(\ell),h(\ell)}(f)$ is defined by~\eqref{eq:swh}.
We introduce the following notation:
\begin{equation}
    k = k(\ell,s) = \begin{cases}
        \ms &\text{if } m(\ell)\geq \ms\\
        m(\ell) &\text{otherwise}
    \end{cases}
\end{equation}
and
\begin{equation}
    \varsigma = \varsigma(\ell,s ) = \begin{cases}
        s &\text{if } m(\ell)\geq \ms\\
        m(\ell)+1 &\text{otherwise}
    \end{cases}
\end{equation}

Remark that, using this notation the kernel $w_{m(\ell)}$ is of order greater than or equal to $k$ and $\varsigma\le s$. Moreover, using classical embedding theorems (see \cite{DINEZZA2012521}), there exists a positive constant $\tilde L$ that depends only on $L$, $s$ and $p$, such that for $\varsigma\in\{2,\ldots,\ms\}$, we have $\tilde\sobolev_{s,p}(L)\subset\tilde\sobolev_{\varsigma,p}(\tilde L)$. For $\varsigma=s$ we also denote $\tilde L=L$.

Now, denoting $h=h(\ell)$ and using a Taylor expansion of $f$, we obtain:
\begin{equation}
    S_{W(\ell),h(\ell)}(f) \leq (k\vee1) \|W(\ell)\|_{\infty}
    \left(\sum_{\vert\alpha\vert=k} I_\alpha\right)^{1/p}
\end{equation}
where
 \begin{align}
     I_\alpha
     &= h^{pk}\int_{\Delta_{d,\mathbf 0}} \int_{\Delta_{d}} \int_0^1
    \left\vert
    \left(D^\alpha f(t+\tau hu)-D^\alpha f(t)\right)
    \right\vert^p \D \tau\D u \D t\\
    &\leq h^{pk}\int_{\Delta_{d,\mathbf 0}} \int_{\Delta_{d}} \int_0^1
    \|hu\|_2^{d+p(\varsigma-k)}\frac{\left\vert D^\alpha f(t+\tau hu)-D^\alpha f(t)\right\vert^p}{\|\tau hu\|_2^{d+p(\varsigma-k)}}
    \D \tau\D u \D t\\
    &\leq d^{(d+p)/2} h^{p\varsigma}
    \int_0^1 \int_{\Delta_{d,\mathbf 0}} \int_{\Delta_{d}}
    \frac{\left\vert D^\alpha f(x)-D^\alpha f(t)\right\vert^p}{\|x-t\|_2^{d+p(\varsigma-k)}} \D x \D t\D\tau\\
     &\leq d^{(d+p)/2} \tilde{L}^p h^{p\varsigma}.
\end{align}
We thus obtain
\begin{equation}\label{eq:100}
    \|\e_f^n \hat f_{\ell}-f\|_p \leq \big[C(d,p,s)\|W(\ell)\|_{\infty}\big]\, \tilde{L} h^\varsigma
\end{equation}
where
\begin{equation}
    C(d,p,s)= (2^dd^{\frac{d+p}2})^{1/p} (\ms\vee 1).
\end{equation}

If $m(\ell)\geq\ms$, since $\tilde L=L$ and $\varsigma=s$, we deduce from \eqref{eq:100} that:
\begin{equation}\label{eq:101}
    \|\e_f^n \hat f_{\ell}-f\|_p \leq \big[C(d,p,s)\|W(\ell)\|_{\infty}\big]\, {L} (h(\ell))^s.
\end{equation}

Assume now that $m(\ell)<\ms$. Then $\varsigma=m(\ell)+1$ and \eqref{eq:100} writes
\begin{equation}
    \|\e_f^n \hat f_{\ell}^{\mathrm{iso}}-f\|_p \leq \big[C(d,p,s) \|W(\ell)\|_{\infty}\big]\, \tilde{L} (h(\ell))^{m(\ell)+1}.
\end{equation}
Remark that
\begin{align}
   (h(\ell))^{m(\ell)+1} &= \exp(-\ell(m(\ell)+1)) \\
   &\leq \exp\left(-\ell\left(\frac{\log n}{2\ell}+\frac12\right)\right)\\
   &\leq \sqrt{\frac{h_n^*}{n}}.
\end{align}
Thus, using Lemma~\ref{lem:norme-omega-m}, for $m(\ell)<\ms$ we obtain:
\begin{equation}
    \|\e_f^n \hat f_{\ell}^{\mathrm{iso}}-f\|_p \leq \big[C(d,p,s) (\ms+1)^{3d/2}\big]\,
    \tilde{L}\sqrt{\frac{h_n^*}{n}} \label{eq:102}.
\end{equation}
Combining \eqref{eq:101} and~\eqref{eq:102} we obtain the proposition.

\subsection{Proof of Proposition~\ref{prop:majorant}}

In the following, $\mathcal L$ is either $\mathcal L_{\mathrm{ani}}$ or $\mathcal L_{\mathrm{iso}}$ and $\hat f_\ell$ then denotes $\hat f_\ell^{\mathrm{ani}}$ or $\hat f_\ell^{\mathrm{iso}}$.
Let $\ell\in\mathcal{L}$. We define

 \begin{equation}
  M_p(\ell) = \frac 1{\sqrt{nV_{\h(\ell)}}} \sum_{\varepsilon\in\{0,1\}^d} \Gamma_\varepsilon(W(\ell), h(\ell), p).
 \end{equation}

 First, assume that $1\leq p\leq 2$.  In this case $M_p(\ell)=\hat M_p(\ell)$, which implies that
 \begin{align}
     \e_f^n \left\{\|{\hat{f}_\ell-\e_f^n\hat{f}_\ell}\|_p-(1+\tau)\hat{M}_p(\ell)\right\}_+^q
     &\leq A_{p,q}(\ell)
 \end{align}
 where
 \begin{equation}
     A_{p,q}(\ell) = \e_f^n \left\{\|{\hat{f}_\ell-\e_f^n \hat{f}_\ell}\|_p- (1+\tau/2)M_p(\ell)\right\}_+^q.
 \end{equation}
 {Next, assume that $p>2$. Consider the event {\begin{equation}
   \label{eq:17}
   \mathcal D_{\ell} = \left\{\sum_{\varepsilon\in\{0,1\}^d}\|{\xi_{W^*(\ell),h(\ell)}}\|_{p/2,\varepsilon}^{1/2}\leq \delta2^d\left(nV_{h(\ell)}\right)^{1/4}\right\}
 \end{equation}}
 with  $\delta=\frac{\tau}{2(1+\tau)}$.} We have
 \begin{align}
    \MoveEqLeft{\left\{\|{\hat{f}_\ell-\e_f^n\hat{f}_\ell}\|_p-(1+\tau)\hat{M}_p(\ell)\right\}_+}\\
    &= \left\{\|{\hat{f}_\ell-\e_f^n\hat{f}_\ell}\|_p-(1+\tau)\hat{M}_p(\ell)\right\}_+\1_{\mathcal{D}_{\ell}} \\
    &\qquad+\left\{\|{\hat{f}_\ell-\e_f^n\hat{f}_\ell}\|_p-(1+\tau/2)M_p(\ell) + (1+\tau/2)M_p(\ell)-(1+\tau)\hat{M}_p(\ell)\right\}_+\1_{\bar{\mathcal{D}}_{\ell}}\\
     &\leq \left\{\|{\hat{f}_\ell-\e_f^n\hat{f}_\ell}\|_p-(1+\tau/2){M_p(\ell)}\right\}_+\\
     &\qquad+(1+\tau/2){M_p(\ell)}\1_{\bar {\mathcal D}_{\ell}}+
     \left\{\|{\hat{f}_\ell-\e_f^n\hat{f}_\ell}\|_p-(1+\tau)\hat{M}_p(\ell)\right\}_+\1_{\mathcal{D}_{\ell}}.
 \end{align}
 Last inequality is true since $\hat{M}_p(\ell)\geq0$. This implies:
 \begin{align}
     \e_f^n \left\{\|{\hat{f}_\ell-\e_f^n\hat{f}_\ell}\|_p-(1+\tau)\hat{M}_p(\ell)\right\}_+^q
     \leq 3^{q-1} \left(A_{p,q}(\ell)+B_{p,q}(\ell)+C_{p,q}(\ell)\right)
 \end{align}
 where
 \begin{equation}
     B_{p,q}(\ell) =  (1+\tau/2)^q(M_p(\ell))^q\p_f^n\left({\bar {\mathcal D}_{\ell}}\right)
 \end{equation}
 and
 \begin{equation}
   \label{eq:23}
   C_{p,q}(\ell)=\e_f^n \left(\left\{\|{\hat{f}_\ell-\e_f^n\hat{f}_\ell}\|_p-(1+\tau)\hat{M}_p(\ell)\right\}_+^q\1_{\mathcal{D}_{\ell}}\right).
 \end{equation}

\paragraph{Control of $\boldsymbol{A_{p,q}(\ell)}$.}

 Remark that
 \begin{align}
     A_{p,q}(\ell) &\le \e_f^n \left\{\sum_{\varepsilon\in\{0,1\}^d}\|{\hat f_{\ell}-\e_f^n \hat f_{\ell}}\|_{p,\varepsilon}-{\frac{(1+\tau/2)\Gamma_\varepsilon(W(\ell), h(\ell),p)}{\sqrt{nV_{h(\ell)}}}}\right\}_+^q\\
     &\leq 2^{d(q-1)} \sum_{\varepsilon\in\{0,1\}^d} \mathcal I_{q,\varepsilon},
 \end{align}
 where
 \begin{equation}
     \mathcal I_{q,\varepsilon} = \e_f^n \left\{\|{\hat{f}_\ell-\e_f^n \hat{f}_\ell}\|_{p,\varepsilon}-{\frac{(1+\tau/2)\Gamma_\varepsilon(W(\ell), h(\ell),p)}{\sqrt{nV_{h(\ell)}}}}\right\}_+^q
 \end{equation}
Thus, using Lemma~\ref{lem:espY} and Lemma~\ref{lem:bousq} with  $r=p$ we can write:
 \begin{align}\label{controlapq}
     \left(nV_{h(\ell)}\right)^{q/2}\mathcal I_{q,\varepsilon} &= \e_f^n \left\{\|\xi_{W(\ell),h(\ell)}\|_{p,\varepsilon}-(1+\tau/2)\Gamma_{\varepsilon}(W(\ell),h(\ell),p)\right\}_+^q\\
     &\leq q\int_0^{+\infty} y^{q-1}\p_f^n \left(\|\xi_{W(\ell),h(\ell)}\|_{p,\varepsilon}-(1+\tau/2)\Gamma_{\varepsilon}(W(\ell),h(\ell),p)>y\right) \D y\\
     &\leq q\int_0^{+\infty} y^{q-1}\p_f^n \left(\|\xi_{W(\ell),h(\ell)}\|_{p,\varepsilon}-\e_f^n \|\xi_{W(\ell),h(\ell)}\|_{p,\varepsilon}> \frac{\tau}{2}\Gamma_{\varepsilon}(W(\ell),h(\ell),p)+y\right) \D y\\
    &{\leq q\exp(-C_1\alpha_n(p))\int_0^{+\infty} y^{q-1}\exp\left(-\frac{C_2 y^2(\alpha_n(p))^{-1}}{\|W(\ell)\|_{2\vee p}^2+y\|W(\ell)\|_p}\right) \D y.}\label{equat}
 \end{align}
Using Lemma~\ref{lem:norme-omega-m}, Condition \eqref{cond:ml} on $m(\ell)$, we have for $y\ge 1$ \begin{equation}
 \frac{C_2 y^2(\alpha_n(p))^{-1}}{\|W(\ell)\|_{2\vee p}^2+y\|W(\ell)\|_p}\ge \frac{C_2 y(\beta_n(p))^{-1}}{2},
 \end{equation}
 where $\beta_n(p)=\alpha_n(p)\left(\log n+3/2\right)^{4d}$.
 Using this bound and doing the change of variable $z=(\beta_n(p))^{-1} y$ in \eqref{equat} we obtain that
 \begin{align}
     \left(nV_{h(\ell)}\right)^{q/2}\mathcal I_{q,\varepsilon} &\le C(\beta_n(p))^q\exp(-C_1\alpha_n(p))
     \end{align}
 where $C$ depends only on $C_2$ $p$ and $q$.

 This implies that
 \begin{equation}
   \label{eq:7}
   A_{p,q}(\ell) = \mathcal O(n^{-q}).
 \end{equation}

 \paragraph{Control of $\boldsymbol{B_{p,q}(\ell)}$.}

Here, we consider $p>2$.\\
Let $\ell\in\mathcal{L}_{\mathrm{ani}}$. Since $h(\ell)$ satisfies $nV_{h(\ell)}\ge (\log n)^c$, using Lemma~\ref{lem:espY}, there exists $N_0=N_0(c,\tau,F_\infty,W^{\circ})$ such that for any $n\geq N_0$:
 \begin{align}
   (1+\tau/2)\Gamma_\varepsilon(W^*(\ell),h(\ell),p/2)\leq \delta^2(nV_{h(\ell)})^{1/2}.
 \end{align}
Let $\ell\in\mathcal{L}_{\mathrm{iso}}$. Using Lemma~\ref{lem:norme-omega-m} and ~\ref{lem:espY}, we have
 \begin{align}
   (1+\tau/2)\Gamma_\varepsilon(W^*(\ell),h(\ell),p/2)&\leq (1+\tau/2)C_0\|W^*(\ell)\|_{2\vee \frac{p}{2}}\\
   &\leq (1+\tau/2)C_0 \left(\frac{\|W(\ell)\|_{4\vee p}}{\|W(\ell)\|_{2}} \right)^2\\
   &  \leq (1+\tau/2)C_0 (m(\ell)+1)^{2d} \\
   & \leq  (1+\tau/2)C_0 \left(\frac{\log n}{2\ell}+\frac{3}{2}\right)^{2d}\\
   &\leq  2^{2d-1} (1+\tau/2)C_0\left[  \left(\frac{\log n}{2\ell}\right)^{2d}+ \left(\frac{3}{2}\right)^{2d}\right]\\
   &\leq  2^{-1} (1+\tau/2)C_0\left[  \sqrt{nV_{h(\ell)}}\left(\frac{ n^{-\frac{1}{4d}}\log n}{\ell e^{-\ell/4}}\right)^{2d}+ 3^{2d}\right]\label{cond:hmin1}
   \end{align}
 Since the bandwidth $h(\ell)$ satisfies $nV_{h(\ell)}\ge (\log n)^c$, this implies that $\ell\le \ell_{max}$ where
 \begin{equation*}
 \ell_{max}= \left[\frac{1}{d}\left(\log n-c\log \log n\right)\right].
 \end{equation*}
  As a consequence
  \begin{align}\frac{ n^{-\frac{1}{4d}}\log n}{\ell e^{-\ell/4}}\le&  \max\left( \frac{ n^{-\frac{1}{4d}}\log n}{ e^{-1/4}}, \frac{ n^{-\frac{1}{4d}}\log n}{\ell_{max} e^{-\ell_{max}/4}} \right)\\
 \le & \max\left( \frac{ n^{-\frac{1}{4d}}\log n}{ e^{-1/4}}, (\log n)^{-\frac{c}{4d}}\frac{\log n}{\ell_{max}}\right)\label{cond:hmin2}
  \end{align}
   Using $nV_{h(\ell)}\ge (\log n)^c$, \eqref{cond:hmin1} and \eqref{cond:hmin2}, we deduce that there exists $N_0=N_0(c,\tau,F_\infty,\delta)$ such that for any $n\geq N_0$:
  \begin{align}
   (1+\tau/2)\Gamma_\varepsilon(W^*(\ell),h(\ell),p/2)&\leq \delta^2(nV_{h(\ell)})^{1/2}.
 \end{align}

As a consequence (in both cases $\ell\in\mathcal{L}_{\mathrm{ani}}$ and $\ell\in\mathcal{L}_{\mathrm{iso}}$), we have for $n\ge N_0$
 \begin{align}\label{eq:Dell}
    \p \left(\overline{\mathcal{D}}_{\ell}\right)
    &\leq \sum_{\varepsilon\in\{0,1\}^d}\p_f^n\left(\|{\xi_{W^*(\ell),h(\ell)}}\|_{{ \frac{p}{2},\varepsilon}}\geq\delta^2(nV_{h(\ell)})^{1/2}\right)\\
    &\leq \sum_{\varepsilon\in\{0,1\}^d}\p_f^n\left(\|{\xi_{W^*(\ell),h(\ell)}}\|_{{ \frac{p}{2},\varepsilon}}\geq(1+\tau/2)\Gamma_\varepsilon(W^*(\ell),h(\ell),p/2)\right)\\
    &\leq  2^d \exp\{-C_1\alpha_n(p/2)\},
 \end{align}
where the last line is a consequence of Lemma \ref{lem:espY} and~\ref{lem:bousq}. Then, using that $nV_{h(\ell)}\ge (\log n)^c$, we have
 \begin{align}
   B_{p,q}(\ell)& \leq C\left( \log n \right)^{-q/c}\exp\{-C_1\alpha_n(p/2)\},
 \end{align}
 for some positive constant $C$ that depends on $\tau$, $F_\infty$, $p$  and $q$ (and $W^{\circ}$ in anisotropic case). This leads finally to
\begin{equation}B_{p,q}(\ell)=\mathcal O(n^{-q}).
\end{equation}

 It remains to upper bound $C_{p,q}(\ell)$ for $q\geq1$ and $p>2$.

 \paragraph{Control of $\boldsymbol{C_{p,q}(\ell)}$.}
 Recall that $p\geq 2$. Let us remark that
 \begin{align}
   \left|\hat{M}_p(\ell) - M_p(\ell)\right|
   &= \left|\sum_{\varepsilon\in\{0,1\}^d} \frac{C_p^* \|W(\ell)\|_{2}}{(nV_{h(\ell)})^{1/2}} Z_\varepsilon(\ell,p)\right|
 \end{align}
 where
 \begin{equation}
 Z_\varepsilon(\ell,p)=\left(\int_{\Delta_{d,\varepsilon}} \left(\e_f^n  \mathcal{K}_{W^*(\ell),h(\ell)}(t,X_1)         \right)^{p/2}    dt  \right)^{1/p}   -\left(\int_{\Delta_{d,\varepsilon}} \left( \frac{1}{n}\sum_{j=1}^n \mathcal{K}_{W^*(\ell),h(\ell)}(t,X_j)          \right)^{p/2}   dt   \right)^{1/p}.
 \end{equation}
 We have
 \begin{align}
 |Z_\varepsilon(\ell,p)|=& \left|\sqrt{ \left\| \e_f^n \mathcal{K}_{W^*(\ell),h(\ell)}(\cdot,X_1)  \right\|_{p/2,\varepsilon}  }-\sqrt{ \left\| \frac{1}{n}\sum_{j=1}^n \mathcal{K}_{W^*(\ell),h(\ell)}(\cdot,X_j)   \right\|_{p/2,\varepsilon}  }\right|\\
  \leq& (nV_{h(\ell)})^{-1/4}\|{\xi_{W^*(\ell),h(\ell)}(\cdot)}\|_{p/2,\varepsilon}^{1/2}.
 \end{align}
 This implies that
 \begin{align}\label{eq:24}
   \left|\hat{M}_p(\ell) - M_p(\ell)\right|
   &\leq  \frac{C_p^*\|W(\ell)\|_{2}}{(nV_{h(\ell)})^{3/4}}
   \sum_{\varepsilon\in\{0,1\}^d}\|{\xi_{W^*(\ell),h(\ell)}(\cdot)}\|_{p/2,\varepsilon}^{1/2}.
 \end{align}

 Thus, under $\mathcal{D}_{\ell}$ defined in \eqref{eq:17} we have
 \begin{align}
   \left|\hat{M}_p(\ell) - M_p(\ell)\right|
   &\leq  \frac{2^{d}C_p^*\|W(\ell)\|_{2}}{(nV_{h(\ell)})^{3/4}}
   \delta (nV_{h(\ell)})^{1/4}\\
   &\leq \delta M_p(\ell), 
 \end{align}
 \begin{align}
   \hat{M}_p(\ell)\geq \left(1-\delta\right)M_p(\ell),
 \end{align}
 and, since $(1-\delta)(1+\tau)=1+\tau/2$:
 \begin{align}
   (1+\tau)\hat{M}_p(\ell)\geq (1+\tau/2)M_p(\ell).
 \end{align}
 This implies that
 \begin{align}
   C_{p,q}(\ell) \leq A_{p,q}(\ell) = \mathcal O(n^{-q}).
 \end{align}

\subsection{Proof of Theorem~\ref{thm:oracle}}

{First, we introduce the following notation: for any $\ell,\ellp\in\mathcal L_{\mathrm{ani}}$, we denote $\ell\preceq\ellp$ if, for any $i=1,\ldots,d$, we have $\ell_i\leq\ellp_i$. Let $\ell\in\mathcal L_{\mathrm{ani}}$ be an arbitrary multi-index. To simplify the notation, we use $\hat{f}_\ell=\hat{f}_\ell^{\mathrm{ani}}$ and $\hat{f}=\hat{f}^{\mathrm{ani}}$.
 \\
Using the definition of $\hat B_p(\ell)$,} we easily obtain:
 \begin{align*}
     \|f-\hat{f}\|_p
     &\le \|f-\hat{f}_\ell\|_p+\|\hat{f}_{\hat{\ell}\wedge\ell}-\hat{f}_\ell\|_p+\|\hat{f}_{\hat{\ell}\wedge\ell}-\hat{f}_{\hat \ell}\|_p\\
     &\le \|f-\hat{f}_\ell\|_p+ \hat B_p(\hat{\ell})+ (1+\tau)\widehat{M}_p(\hat{\ell},\ell)+ \hat B_p(\ell)+ (1+\tau)\widehat{M}_p(\ell,\hat{\ell}).
 \end{align*}
 Using the definition of $\hat\ell$, we deduce:
 \begin{align*}
   \|f-\hat{f}\|_p
   &\leq \|f-\hat{f}_\ell\|_p+  2\left(\hat B_p({\ell})+ (1+\tau)\widehat{M}_p(\ell)\right)+2(1+\tau)\widehat{M}_p(\ell\wedge\hat\ell)\\
   &\leq \|f-\hat{f}_\ell\|_p+  2\hat B_p({\ell})+ 4(1+\tau)\max_{\ellp\preceq\ell} M_p(\ellp) + 4(1+\tau)\max_{\ellp\preceq\ell} \left(\hat M_p(\ellp)-M_p(\ellp)\right).
 \end{align*}
 This implies that:
 \begin{align}
     R_n^{(p,q)}(\hat f,f)
     &\leq  R_n^{(p,q)}(\hat f_\ell,f)
     + 2  \left(\efn\hat B_p^q(\ell)\right)^{1/q}\\
     &\qquad+ 4(1+\tau)\left(\efn\max_{\ellp\preceq\ell}|\widehat M_p(\ellp)-M_p(\ellp)|^q\right)^{1/q} \\
     &\qquad+ 4(1+\tau) \max_{\ell'\preceq\ell}{M}_p(\ell').
 \end{align}
 It remains to bound each term of the right hand side of this inequality.

 \paragraph{1.} Remark that, using the triangular inequality, we have:
 \begin{align*}
 \hat B_p(\ell)
     &\le 2 \max_{\ellp\in\mathcal L} \left\{ \|\hat{f}_{\ellp}-\e_f^n \hat{f}_{\ellp}\|_p -(1+\tau)\widehat{M}_p(\ellp)\right\}_+ +\max_{\ellp\in\mathcal L} \| \e_f^n \hat{f}_{\ellp}-\e_f^n \hat{f}_{\ell\wedge\ellp} \|_p.
 \end{align*}
 {This readily implies
 \begin{align}
     \left(\efn\hat B_p^q(\ell)\right)^{1/q}
     &\leq 2\sum_{\ellp\in\mathcal L} \left(\efn\left\{ \|\hat{f}_{\ellp}-\e_f^n \hat{f}_{\ellp} \|_p -(1+\tau)\widehat{M}_p(\ellp)\right\}_+^q\right)^{1/q}\\
     &\qquad + \max_{\ellp\in\mathcal L} \| \e_f^n \hat{f}_{\ellp}-\e_f^n \hat{f}_{\ell\wedge\ellp} \|_p\\
     &\leq 2\const5^{1/q}(\#\mathcal L)n^{-1} + \max_{\ellp\in\mathcal L} \| \e_f^n \hat{f}_{\ellp}-\e_f^n \hat{f}_{\ell\wedge\ellp} \|_p,\label{eq:par1}
 \end{align}}
 where the last inequality follows immediately from Proposition~\ref{prop:majorant}.

 \paragraph{2.}
 For $p\le 2$, we have  $\hat M_p(\ell)-M_p(\ell)=0$.

 Let $p>2$. Here and in the following paragraph, $C$ stands for a constant that depends on $p$, $q$, $\tau$, $F_\infty$ and $W^\circ$ and that can change of values from line to line.
 Using \eqref{eq:24}, we obtain that for $\ellp\preceq\ell$
 \begin{align}
 |\hat M_p(\ellp)-M_p(\ellp)|
  &\leq \frac{C_p^*\|W^\circ\|_{2}}{(nV_{h(\ell)})^{1/2}(nV_{h(\ellp)})^{1/4}}
   \sum_{\varepsilon\in\{0,1\}^d}\|{\xi_{W^*(\ellp),\h(\ellp)}(\cdot)}\|_{p/2,\varepsilon}^{1/2}.
\end{align}
{We have
\begin{align*}
\efn\max_{\ellp\preceq\ell}|\widehat M_p(\ellp)-M_p(\ellp)|^q\le & \efn\left\{\max_{\ellp\preceq\ell}|\widehat M_p(\ellp)-M_p(\ellp)|^qI_{\cap_{\ellp\preceq \ell}\mathcal{D}_{\ellp}}\right\}\\
&+\sum_{\ellp\preceq \ell} \efn\left\{\max_{\ellp\preceq\ell}|\widehat M_p(\ellp)-M_p(\ellp)|^qI_{\overline{\mathcal{D}}_{\ellp}}\right\},
\end{align*}
where the events $\mathcal{D}_\ellp$ are defined by \eqref{eq:17}. Then, using \eqref{eq:Dell} and that $h(\ell)\in\mathcal{H}_n$, we obtain for $n$ large enough that
\begin{align*}
\efn\max_{\ellp\preceq\ell}|\widehat M_p(\ellp)-M_p(\ellp)|^q\le & \left( \frac{C_p^*\|W^\circ\|_{2}\delta 2^d}{(nV_{h(\ell)})^{1/2}}  \right)^q+\sum_{\ellp\preceq \ell} \left( \frac{C_p^*\|W^\circ\|_{\infty} 2^{d+1}}{\sqrt{V_{h(\ell)}}(nV_{h(\ell)})^{1/2}}  \right)^q \p \left(\overline{\mathcal{D}}_{\ellp}\right)\\
\le & \left( \frac{C_p^*\|W^\circ\|_{2}\delta 2^d}{(nV_{h(\ell)})^{1/2}}  \right)^q+\sum_{\ellp\preceq \ell} \left( \frac{C_p^*\|W^\circ\|_{\infty} 2^{d+1}\sqrt{n}}{(\log n)^{2d+1}}  \right)^q 2^d\exp\{-C_1\alpha_n(p/2)\}.
\end{align*}
}
 Now since $\#\mathcal L$ is bounded by $(\log n)^d$, we have
 \begin{align*}
 \left(\efn\max_{\ellp\preceq\ell}|\widehat M_p(\ellp)-M_p(\ellp)|^q\right)^{1/q}&\leq \frac{C}{(nV_{h(\ell)})^{1/2}}.
 \end{align*}

 \paragraph{3.}

By using Lemma~\ref{lem:espY} we obtain
 \begin{align}
     M_p(\ell) = &\frac{1}{\sqrt{nV_{\hell}}} \sum_{\varepsilon\in\{0,1\}^d} \Gamma_\varepsilon(W^\circ, h(\ell), p)\\
    \le  & \frac{C\ \|W^\circ\|_{p\vee 2}}{\sqrt{nV_{\hell}}}.\label{eq:par3}
 \end{align}

This implies that for $\ellp\preceq\ell$
 \begin{equation}
 4 \max_{\ell'\preceq\ell}{M}_p(\ell')\le \frac{C}{\sqrt{nV_{\hell}}}.
 \end{equation}

\subsection{Proof of Theorem~\ref{thm:adaptive-anisotropic}}
The proof of this result is split into two main parts: the proof of the upper bound and the proof of the lower bound.

\paragraph{Upper bound}
Set $s\in\prod_{i=1}^d(0,M_i+1]$. Define $\ell^*(s)=(\ell_1^*(s),\ldots,\ell_d^*(s))$ by:
\begin{equation}
    \ell_i^*(s) = \left\lceil \frac{\bar s}{s_i(2\bar s+1)}\log n\right\rceil,\qquad i=1\ldots,d
\end{equation}
where $\lceil x\rceil$ denotes the least integer greater than or equal to $x$. Note that $h_i(\ell^*)$ is such that
\begin{equation}\label{eq:200}
    \frac{h_i^*(s)}{e} \leq h_i(\ell^*) \leq h_i^*(s)
\end{equation}
where
\begin{equation}
    h_i^*(s) = n^{-\frac{\bar s}{s_i(2\bar s+1)}}.
\end{equation}
This implies that there exists $n_0=n_0(s,p)\in\N$ such that for any $n\geq n_0$ we have $\ell^*\in\mathcal L_{\mathrm{ani}}$.

Combining \eqref{eq:200} with Proposition~\ref{prop:bias-anisotropic} and Theorem~\ref{thm:oracle}, the upper bound follows.

\paragraph{Lower bound}

For the sake of simplicity we only prove the lower bound for anisotropic Sobolev-Slobodetskii classes. Let $(s_1,\dots,s_d)$ be a vector of positive real numbers and let $L>0$. We also assume that $s_i\notin\N$. We intend to prove that $n^{-{\bar s}/{(2\bar s+1)}}$ is the minimax rate of convergence over the class $\sobolev_{s,p}(L)$. To do so, using Lemma~3 in \citet{MR3346701}, we have to construct a family of density functions $\{f_w : w\in\mathcal W_0\}$, indexed by $\mathcal W_0\ni0$, that satisfies the following properties:
\begin{align}
  (a) & \qquad f_w \in \sobolev_{s,p}(L), \quad w\in\mathcal W_0\\
  (b) & \qquad \|f_w-f_{w'}\|_p \geq 2\rho_n, \quad w,w'\in\mathcal W_0\\
  (c) & \qquad  \mathfrak{I} =
    \frac{1}{|\mathcal W|^2} \sum_{w\in\mathcal W} \e_0\left(
    \prod_{k=1}^n\frac{f_w}{f_0}(X_k)
    \right)^2 \leq \mathfrak a,
\end{align}
where $|\mathcal W|$ denotes the cardinality of $\mathcal W=\mathcal W_0\setminus \{0\}$. Under these assumptions  we have:
\[
  \liminf_{n\to\infty}\inf_{\tilde f} \rho_n^{-1}\sup_{f\in\sobolev_{s,p}(L)} \left(\e_f^n\|\tilde f-f\|_p^q\right)^{1/q}
  \geq (\sqrt{\mathfrak a}+\sqrt{\mathfrak a+1})
\]
where the infimum is taken over all possible estimators. This implies the result. It remains to construct such a family. The rest of the proof is decomposed into several steps.

\paragraph{Step 1.} Here, we construct a finite set of functions used in the rest of the proof. We consider two auxiliary functions $\psi:\R\to\R$ and $H:\R\to\R$ defined, for any $u\in\R$ by
\[
	\psi(u)=\exp(-1/(1-u^2))\1_{(-1,1)}(u)
	\quad\text{and}\quad
	H(u) = -\1_{(-1, 0)} + \1_{(0,1)}.
\]
Using these functions, we define, for any $u\in\R$, $\varphi(u)=H\star\psi(2u)$.

For any $i=1,\dotsc,d$, we consider the bandwidth
\[
 	h_i = n^{-\frac{\bar s}{2\bar s+1}\frac{1}{s_i}}
\]
and we set $R_i=1/(2h_i)$. We assume, without loss of generality, that $R_i$ is an integer. Let $\mathcal{R}=\prod_{i=1}^d \{0,\dotsc,R_i-1\}$ and define, for any $r=(r_1,\dotsc,r_d)\in\mathcal{R}$, the function $\phi_r:\Delta_d\to\R$ by:
\[
	\phi_r(x) = \prod_{i=1}^d \varphi\left(\frac{x_i-x_i^{(r)}}{h_i}\right).
\]
where $x_i^{(r)}=(2r_i+1)h_i$. Finally, for any $w:\mathcal{R}\to\{0,1\}$ we define:
\[
	f_w = \1_{\Delta_d} + \rho_n\sum_{r\in\mathcal{R}} w(r) \phi_r
\]
where
\[
	\rho_n
	= c_1
	n^{-\frac{\bar s}{2\bar s+1}}
  \qquad\text{with}\qquad
  c_1= \frac{L}{d\ell(s,p)\Phi^d}
\]
and
\[
  \Phi = \max_{0\leq k\leq \max_i \lfloor s_i\rfloor+1} \|\varphi^{(k)}\|_\infty
  \qquad\text{and}\qquad
  \ell(s,p) = \left(\frac{6\cdot 2^p}{p(1-\sigma)}+8\sum_{k\geq1}(2k)^{-(1+p\sigma)}\right)^{1/p}
\]
where $\sigma=\min\{s_i-\lfloor s_i\rfloor : i=1,\dotsc,d\}$.

\paragraph{Step 2.} We intend to prove that, for $n$ large enough, $f_w$ is a probability density that belongs to  $\sobolev_{s,p}(L)$.

\noindent
\textit{(i)} Remark that, for any $u\in\R$, since $\|H\|_\infty\leq 1$, we have
\begin{align}
	|\varphi(u)|
	&= \left|\int_\R H(2u-v)\psi(v) dv\right| \\
	&\leq \int_{-1}^1 \psi(v) dv \\
	&\leq 2/e.
\end{align}
This implies that $\|\phi_r\|_\infty\leq (2/e)^d$ for any $r\in\mathcal{R}$. Moreover, note that:
\[
	\operatorname{Supp}(\phi_r) = \prod_{i=1}^d [x_i^{(r)}-h_i, x_i^{(r)}+h_i ].
\]
Thus, the Lebesgue measure of $\operatorname{Supp}(\phi_r)\cap\operatorname{Supp}(\phi_{r'})$ is null for $r\neq r'$. This implies that, for $n$ large enough:
\[
	f_w(x)\geq 1 - \rho_n (2/e)^{d}>0,
	\qquad x\in\Delta_d.
\]
\noindent
\textit{(ii)} Remark that $\varphi$ is an odd function such that $\operatorname{Supp}(\varphi)=[-1,1]$. This implies that $\int_\R \varphi(u) du =0$ which also implies, using Fubini's Theorem, that $\int_{\Delta_d} \phi_r(x) dx=0$ for any $r\in\mathcal{R}$. As a consequence we have $\int_{\Delta_d} f_w(x) dx=1$.

Combining points \textit{(i)}\/ and \textit{(ii)}\/ we deduce that $f_w$ is a probability density for any $w:\mathcal{R}\to\{0,1\}$.
It remains to prove that $f_w$ belongs to the anisotropic class $\sobolev_{s,p}(L)$.

\noindent
\textit{(iii)} Set  $i\in\{1,\dotsc,d\}$ and consider $x,y\in\Delta_d$ such that $x_j=y_j$ for any $j\neq i$. For the sake of readability we denote $s_i=m+\alpha$ with $m=\lfloor s_i\rfloor\in\N$ and $0<\alpha<1$. Remark that $\varphi$ is an infinitely differentiable function. Thus:
\begin{align}
	D_i^m f_w(x)-D_i^m f_w(y)
	&= \rho_n \sum_{r\in\mathcal{R}} \prod_{j\neq i}\varphi\left(\frac{x_j-x_j^{(r)}}{h_j}\right)D_i^m\left(\varphi\left(\frac{x_i-x_i^{(r)}}{h_i}\right)-\varphi\left(\frac{y_i-x_i^{(r)}}{h_i}\right)\right)\\
	&= \frac{\rho_n}{h_i^m}\sum_{r\in\mathcal{R}} \prod_{j\neq i}\varphi\left(\frac{x_j-x_j^{(r)}}{h_j}\right)\left(\varphi^{(m)}\left(\frac{x_i-x_i^{(r)}}{h_i}\right)-\varphi^{(m)}\left(\frac{y_i-x_i^{(r)}}{h_i}\right)\right).
\end{align}
This implies that:
\begin{align}
	|D_i^m f_w(x)-D_i^m f_w(y)|
	&\leq \frac{\|\varphi\|_{\infty}^{d-1}\rho_n}{h_i^m}\sum_{r\in\mathcal{R}} \left|\varphi^{(m)}\left(\frac{x_i-x_i^{(r)}}{h_i}\right)-\varphi^{(m)}\left(\frac{y_i-x_i^{(r)}}{h_i}\right)\right|.
\end{align}
Denote $A_s=[2sh_i, (2s+2)h_i]$. We have
\begin{align}
    I_i(D_i^m f_w)
    &\leq
    \frac{\|\varphi\|_{\infty}^{d-1}\rho_n}{h_i^m}
    \left(
    \sum_{s,t=0}^{R_i-1}
    \int_{A_s}\int_{A_t}
    \frac{%
    \left|
    \sum_{r=0}^{R_i-1}
    \left(
    \varphi^{(m)}\left(\frac{x-(2r+1)h_i}{h_i}\right) -
    \varphi^{(m)}\left(\frac{y-(2r+1)h_i}{h_i}\right)
    \right)
    \right|^p
    }%
    {|x-y|^{1+p\alpha}}
    dx dy
    \right)^{1/p}\\
	&\leq
	\frac{\|\varphi\|_{\infty}^{d-1}\rho_n}{h_i^m}
    \left(
    	\Upsilon+\tilde\Upsilon
    \right)^{1/p} \label{eq:thm5-1}
\end{align}
where
\[
	\Upsilon = \sum_{s=0}^{R_i-1}
    \int_{A_s}\int_{A_s}
    \frac{%
    \left|
    \left(
    \varphi^{(m)}\left(\frac{x-(2s+1)h_i}{h_i}\right) -
    \varphi^{(m)}\left(\frac{y-(2s+1)h_i}{h_i}\right)
    \right)
    \right|^p
    }%
    {|x-y|^{1+p\alpha}}
    dx dy
\]
and
\[
	\tilde \Upsilon =
	2 \sum_{s=0}^{R_i-1}
	\sum_{\substack{t=0\\ |t-s|\geq 1}}^{R_i-1}
	\int_{A_s}\int_{A_t}
	\frac{%
    \left|
    \varphi^{(m)}\left(\frac{x-(2s+1)h_i}{h_i}\right)
    \right|^p
    }%
    {|x-y|^{1+p\alpha}}
    dy dx
\]

First we control $\Upsilon$.
\begin{align}
    \Upsilon
    &\leq
   	\|\varphi^{(m+1)}\|_\infty^p
    \sum_{s=0}^{R_i-1}
    \int_{A_s}\int_{A_s}
    \frac{%
    \left|\frac{x-y}{h_i}\right|^p
    }%
    {|x-y|^{1+p\alpha}}
    dx dy \\
    &\leq
     \frac{\|\varphi^{(m+1)}\|_\infty^p}%
    {h_i^{p\alpha}}
    \sum_{s=0}^{R_i-1}\int_{A_s}\int_{-2}^2 |u|^{p(1-\alpha)-1} du dv\\
    &\leq
     \frac{\|\varphi^{(m+1)}\|_\infty^p}%
    {h_i^{p\alpha}}
    \int_{0}^1\int_{-2}^2 |u|^{p(1-\alpha)-1} du dv\\
    &\leq
    \frac{2^{1+p(1-\alpha)}\|\varphi^{(m+1)}\|_\infty^p}{p(1-\alpha)}
    h_i^{-p\alpha}.\label{eq:thm5-2}
\end{align}

Now, we control $\tilde\Upsilon$. Note that the sum over $t$ can be decomposed into two different terms for $|t-s|=1$ or $|t-s|\geq 2$.

First, remark that if $x\in A_s$, $y\in A_t$ and $|s-t|\geq 2$ then $|x-y|\geq 2(|s-t|-1)h_i$. This implies that:
\begin{align}
    \tilde\Upsilon_2
    &= 2 \sum_{s=0}^{R_i-1}
	\sum_{\substack{t=0\\ |t-s|\geq 2}}^{R_i-1}
	\int_{A_s}\int_{A_t}
	\frac{%
    \left|
    \varphi^{(m)}\left(\frac{x-(2s+1)h_i}{h_i}\right)
    \right|^p
    }%
    {|x-y|^{1+p\alpha}}
    dy dx \\
	&\leq
	4 \sum_{s=0}^{R_i-1}
	\sum_{k\geq 1}
	\frac{
    (2h_i)^2\|\varphi^{(m)}\|_\infty^p
    }%
    {(2kh_i)^{1+p\alpha}}\\
	&\leq
	\frac{8\|\varphi^{(m)}\|_\infty^p}{h_i^{p\alpha}}
	\sum_{k\geq 1}(2k)^{-(1+p\alpha)}.\label{eq:thm5-3}
\end{align}

Now, it remains to consider the case $|s-t|=1$.Assume first that $t=s+1$. We consider the point $z=(2s+2)h_i$ that satisfies: $|x-z|\leq|x-y|$, $|y-z|\leq|x-y|$ and $z\in A_s\cap A_{s+1}$. We can also remark that $\varphi^{(m)}((z-(2s+1)h_i)/h_i)=0$. We this in mind remark that:
\begin{align}
    \tilde \Upsilon_1^{+}
    &=
	2 \sum_{s=0}^{R_i-2}
	\int_{A_s}\int_{A_{s+1}}
	\frac{%
    \left|
    \varphi^{(m)}\left(\frac{x-(2s+1)h_i}{h_i}\right)
    \right|^p
    }%
    {|x-y|^{1+p\alpha}}
    dy dx\\
	&\leq
   	2\|\varphi^{(m+1)}\|_\infty^p
    \int_{A_s}
    \frac{%
    \left|\frac{x-z}{h_i}\right|^p
    }%
    {|x-z|^{1+p\alpha}}
    dx dy \\
	&\leq
    \frac{2^{1+p(1-\alpha)}\|\varphi^{(m+1)}\|_\infty^p}{p(1-\alpha)}
    h_i^{-p\alpha}.	\label{eq:thm5-4}
\end{align}
In the same way for $t=s-1$:
\begin{align}
    \tilde \Upsilon_1^{-}
    &=
	2 \sum_{s=1}^{R_i-1}
	\int_{A_s}\int_{A_{s-1}}
	\frac{%
    \left|
    \varphi^{(m)}\left(\frac{x-(2s+1)h_i}{h_i}\right)
    \right|^p
    }%
    {|x-y|^{1+p\alpha}}
    dy dx\\
	&\leq
    \frac{2^{1+p(1-\alpha)}\|\varphi^{(m+1)}\|_\infty^p}{p(1-\alpha)}
    h_i^{-p\alpha}.	\label{eq:thm5-5}
\end{align}
Using that $\tilde\Upsilon=\tilde\Upsilon_1^++\tilde\Upsilon_1^-+\tilde\Upsilon_2$ combined with \eqref{eq:thm5-1}, \eqref{eq:thm5-2}, \eqref{eq:thm5-3}, \eqref{eq:thm5-4} and \eqref{eq:thm5-5} leads to:
\[
  \sum_{i=1}^d I_i(D_i^{(\lfloor s_i\rfloor) f_w})
  \leq d\ell(s,p)\Phi^d c_1
  \leq L,
\]
which implies that $f_w$ is a probability density that belongs to $\sobolev_{s,p}(L)$.

\paragraph{Step 3.} To define the set $\mathcal W_0$ we introduce the following notations. Let \[
	c_2=\min \left\{\frac{2^{-d}}{2+4\exp\left(2c_1^2\|\varphi\|_\infty^{2d}\right)}, 2^{p+1}\|\varphi\|_p^d, 2^{-d}/10 \right\}
\]
and define $M=\prod_{i=1}^d R_i = (2^dV_h)^{-1}$ and $m=c_2V_h^{-1}$. Without loss of generality we assume that both $M$ and $m\geq 4$ are integers.
Using Lemma~A3 in \citet{rigollet-tsybakov-2011}, there exists $\mathcal{W}\subset\{w\colon \mathcal{R}\to\{0,1\}\}$ such that:
\begin{itemize}
	\item We have $|\mathcal{W}|\geq 2^{-m}(M/m-1)^{m/2}$.
	\item For any $w\in\mathcal W$ we have:
	\[
		|w| = \sum_{r} w(r) = m
	\]
	\item For any $w,w'\in\mathcal W$, we have:
	\[
		\sum_{r\in\mathcal R} |w(r)-w'(r)| \geq m/2.
	\]
\end{itemize}
Then, define $\mathcal W_0 = \mathcal W\cup\{0\}$. Remark that the last point remains valid if one replaces $\mathcal{W}$ by $\mathcal{W}_0$ thanks to the second point. Remark also that $f_0 \equiv \1_{\Delta_d}$.

\noindent
\textit{(i)} Let $w$ and $w'$ in $\mathcal W_0$.
\begin{align}
    \|f_w - f_{w'}\|_p
    &= \rho_n \left\|\sum_{r\in\mathcal{R}} \big(w(r)-w'(r)\big) \phi_r\right\|_p\\
    &= \rho_n \left(\sum_{s\in\mathcal{R}}\int_{\operatorname{Supp}(\phi_s)}
    \left|\sum_{r\in\mathcal{R}} \big(w(r)-w'(r)\big) \phi_r(u)\right|^p du
    \right)^{1/p}
\end{align}
Using the fact that the functions $\phi_r$ and $\phi_s$ have disjoint supports for $r\neq s$, we have:
\begin{align}
    \|f_w - f_{w'}\|_p
    &=  \rho_n \left(\sum_{s\in\mathcal{R}} \big|w(s)-w'(s)\big|^p \cdot \|\phi_s\|_p^p\right)^{1/p}\\
    &= \rho_n \left(\sum_{r\in\mathcal{R}} \big|w(r)-w'(r)\big|\right)^{1/p}V_h^{1/p}\|\varphi\|_p^d\\
    &\geq \rho_n (c_2/2)^{1/p} \|\varphi\|_p^d\\
    &\geq 2\rho_n.
\end{align}

\noindent\textit{(ii)}
In what follows, we denote by $\e_0$ the expectation under the uniform distribution on $\Delta_d$, with density $f_0$.
\begin{align}
    \mathfrak{I}
    &= \frac{1}{|\mathcal W|^2} \sum_{w\in\mathcal W}
    \left(\int_{\Delta_d} f_w^2(x) dx\right)^n \\
    &\leq \frac{1}{|\mathcal W|^2} \sum_{w\in\mathcal W}
    \left(\int_{\Delta_d} \bigg(1+2\rho_n\sum_{r\in\mathcal{R}} w(r)\phi_r(x) + \rho_n^2\sum_{r,r'\in\mathcal{R}} w(r)w(r')\phi_r(x)\phi_{r'}(x)\bigg) dx\right)^n \\\\
    &\leq \frac{1}{|\mathcal W|}  \big(1+\rho_n^2 mV_h \|\varphi\|_2^{2d}\big)^n.
\end{align}
Last inequality comes from the facts that $\int_{\Delta_d}\phi_r(x) dx=0$ and that the set $\operatorname{Supp}(\phi_r)\cap\operatorname{Supp}(\phi_{r'})$ is negligible (in terms of Lebesgue measure) for $r\neq r'$.
We thus obtain:
\begin{align}
    \mathfrak{I}
    &\leq  \exp\big(n\rho_n^2 c_2 \|\varphi\|_2^{2d}-\log(|\mathcal W|\big)\\
    &\leq \exp\left(c_1^2c_2\|\varphi\|_2^{2d} n^{1/(2\bar s+1)}-\frac{m}{2}\log\left(\frac14\left(\frac{M}{m}-1\right)\right)\right)\\
    &\leq \exp\left(
      \frac{m}{2}
      \left(
        2c_1^2\|\varphi\|_2^{2d}-\log\frac{2^{-d}-c_2}{4c_2}
      \right)
    \right).
\end{align}
Using the definition of $c_2$ we remark that the exponent is nonpositive. This implies that $\mathfrak{J}\leq 1$. Taking all together, the assumptions of Lemma~3 in \citet{MR3346701} are satisfied. Theorem is then proved.

\subsection{Proof of Theorem~\ref{thm:oracle2}}

Let $\ell\in\mathcal{L}_{\mathrm{iso}}$.
{We have
\begin{equation*}
 \|f-\hat{f}^{\mathrm{iso}}\|_p
   \leq \|f-\hat{f}_\ell^{\mathrm{iso}}\|_p+\|\hat{f}_{\hat{\ell}\wedge\ell}^{\mathrm{iso}}-\hat{f}_\ell^{\mathrm{iso}}\|_p+\|\hat{f}_{\hat{\ell}\wedge\ell}^{\mathrm{iso}}-\hat{f}_{\hat \ell}^{\mathrm{iso}}\|_p,
 \end{equation*}
 Note that if $\ell\geq \hat \ell$, then
 \begin{align*}
   \|f-\hat{f}^{\mathrm{iso}}\|_p
   &\leq \|f-\hat{f}_\ell^{\mathrm{iso}}\|_p+\|\hat{f}_{\hat{\ell}\wedge\ell}^{\mathrm{iso}}-\hat{f}_\ell^{\mathrm{iso}}\|_p\\
   &\leq \|f-\hat{f}_\ell^{\mathrm{iso}}\|_p+ \hat B_p ({\hat\ell})+ (1+\tau)\hat{M}_p (\hat\ell,\ell)\\
   &\leq \|f-\hat{f}_\ell^{\mathrm{iso}}\|_p+ \hat B_p ({\hat\ell})+ (1+\tau)\hat{M}_p (\ell)
     + (1+\tau)\hat{M}_p (\hat\ell)\\
   &\leq \|f-\hat{f}_\ell^{\mathrm{iso}}\|_p+ 2\left( \hat B_p (\ell) + (1+\tau)\hat{M}_p (\ell)\right).
 \end{align*}
 Last inequality comes from the definition of $\hat\ell$.
 It is easily seen that the same bound remains valid if  $\ell\leq \hat \ell$. This implies that
 \begin{align}
     R_n^{(p,q)}(\hat{f}^{\mathrm{iso}},f)
     &\leq  R_n^{(p,q)}(\hat{f}_\ell^{\mathrm{iso}},f)
     + 2  \left(\efn\hat B_p^q(\ell)\right)^{1/q}\\
     &\qquad+ 2(1+\tau)\left(\efn|\widehat M_p(\ell)-M_p(\ell)|^q\right)^{1/q} \\
     &\qquad+ 2(1+\tau) {M}_p(\ell).
 \end{align}
}

Following the same arguments of the proof of Theorem~\ref{thm:oracle} (see the second paragraph), we have
\begin{align}
  \left(\e |\hat M_p(\ell)-M_p(\ell)|^q\right)^{1/q}&\leq C \frac{\|W(\ell)\|_{2}}{(nV_{h(\ell)})^{1/2}}.\label{eq:par2}
\end{align}
Applying \eqref{eq:par1}, \eqref{eq:par3} and \eqref{eq:par2}, we deduce the oracle inequality of Theorem~\ref{thm:oracle2}.

\subsection{Proof of Theorem~\ref{thm:adaptive-isotropic}}

First, we prove the upper bound. Set $s>0$. Define:
\begin{equation}
    \ell^*(s) = \left[ \frac{1}{2s+d} \log n \right]
\qquad
\text{and}
\qquad
    h_n^*(s) = n^{-\frac{1}{2s+d}}.
\end{equation}
We note that there exists $n_1=n_1(s,p)$ such that for any $n\geq n_1$ we have $\ell^*(s)\in\mathcal L_{\mathrm{iso}}$ and $1\le \frac{\log n}{2(2s+d)}$. Then we have
\begin{equation}\label{eq:202}
    s\leq m(\ell^*(s))\leq 2s+d+\frac{1}{2}
\end{equation}
and
\begin{equation}\label{eq:201}
    h_n^*(s) \leq h(\ell^*(s)) \leq e h_n^*(s).
\end{equation}

Then using Lemma~\ref{lem:norme-omega-m}, \eqref{eq:202} implies
\begin{equation}\label{eq:203}
    \max(\|W(\ell^*(s))\|_{2\vee p}, \|W(\ell^*(s))\|_\infty)
    \leq \left(2s+d+\frac{3}{2}\right)^{2d}.
\end{equation}

Using~\eqref{eq:201} and ~\eqref{eq:203} in combination with Proposition~\ref{prop:bias-isotropic} and Theorem~\ref{thm:oracle2} entail to the upper bound.
To prove the lower bound, the methodology and construction proposed in the proof of Theorem~\ref{thm:adaptive-anisotropic} are unchanged (we just consider $s_i=s$ for any $i=1,\dotsc,d$). However it remains to prove that the functions $f_w$ defined previoulsely belong to the isotropic Sobolev-Slobodetskii class $\tilde\sobolev_{s,p}(L)$. This is left to the reader.

\subsection{Proof of Lemma~\ref{lem:norme-omega-m}}

Let $m\in\N$. Denote $z_m(u)=\sum_{r=0}^m a_r^{(m)}u^r$.
The solution of the minimization problem~\eqref{eq:2} can be found explicitly and the Lagrangian condition implies the solution is $z_m$. Now, for $p=2$, remark that:
 \begin{equation}
     \|{z_m}\|_2^2 =  (a^{(m)})^\top H_{m} a^{(m)} = (e_0^{(m)})^\top H_{m}^{-1} e_0^{(m)} = (m+1)^2.
 \end{equation}
 Now we will prove that $w_m=z_m$.
The polynomial $z_m$ can be decomposed in the basis $\{\varphi_r, r=0,\ldots,m\}$
as
$$z_m(u)=\sum_{r=0}^m b_r\varphi_r(u).$$
Since $z_m$ is of order $m$, we have
$$b_r=\int_0^1 z_m(u) \varphi_r(u)du=\varphi_r(0)$$
which implies that $w_m=z_m$

Finally we have for $u\in\Delta_1$
 \begin{align}\label{eq:upsilon-1}
     |{w_m}(u)| & \leq \sum_{r=0}^m \sqrt{2r+1}|Q_r(-1)|\sqrt{2r+1}|Q_r(2u-1)|\\
     &\leq \sum_{r=0}^m 2r+1\\
     &= (m+1)^2
 \end{align}
 since $|Q_r(u)|\le |Q_r(-1)|=1$.
 Moreover ${w_m}(0)=\sum_{r=0}^m (2r+1)(Q_r(-1))^2=(m+1)^2$ which implies that $\|w_m\|_\infty=(m+1)^2$.

 \subsection{Proof of Lemma~\ref{lem:espY}}
{ Using Jensen inequality, we have
\begin{align*}
\e_f^n \|\xi_{W,h}\|_{2,\varepsilon}&\leq \left(\frac{   V_h}{n}\right)^{1/2} \left(\int_{\Delta_{d,\varepsilon}}\e_f^n
    \left( \sum_{j=1}^n \mathcal{K}_{W,h}(t,X_j) -\e_f^n
    \mathcal{K}_{W,h}(t,X_j) \right)^{2}dt\right)^{1/2}\\
  &\leq  \sqrt{V_h}  \left(\int_{\Delta_{d,\varepsilon}}
    \e_f^n  \mathcal{K}^2_{W,h}(t,X_1)dt\right)^{1/2}\\
  & \le \left(\int_{\Delta_{d,\varepsilon}}
   \int_{\Delta_{d}}\frac{1}{V_h} \prod_{i=1}^d W_i^2\left(\sigma(t_i)\frac{t_i-x_i}{h_i}\right)f(x)dx
   dt\right)^{1/2}
\end{align*}
Then using a change of variables, we deduce
\begin{equation}
\e_f^n \|\xi_{W,h}\|_{2,\varepsilon} \le \left(\|W\|_2^2\int_{\Delta_{d}} f(x)dx\right)^{1/2}\le \|W\|_2.
\end{equation}
 For $r\le 2$, since the Lebesgue measure of $\Delta_{d,\varepsilon}$ equals to $2^{-d}$, we have using H\"older inequality
 \begin{equation}
   \label{eq:8}
   \e_f^n \|\xi_{W,h}\|_{r,\varepsilon}\leq 2^{-\frac{d(2-r)}{2r}}\e_f^n \|\xi_{W,h}\|_{2,\varepsilon}\leq2^{-\frac{d(2-r)}{2r}}\|W\|_2.
 \end{equation}}
 Let us now assume that $r>2$. Using the Rosenthal inequality we have
 \begin{align}
     \e_f^n |\xi_{W,h}(t)|^r
     &\leq (C_r^*)^r(V_h)^{r/2} \left\{
     \left( \e_f^n \mathcal K_{W,h}^2(t,X_1)\right)^{r/2}
     + 2^{r+1} n^{1-r/2} \e_f^n |\mathcal K_{W,h}(t,X_1)|^r
     \right\}.
 \end{align}
 Using Jensen and Young inequalities we obtain:
 \begin{align}
    \e_f^n \|\xi_{W,h}\|_{r,\varepsilon}
     &\leq  \left(\int_{\Delta_{d,\varepsilon}} \e_f^n |\xi_{W,h}(t)|^r \D t\right)^{1/r}\\
     &\leq C_r^* \left\{\Lambda_{\varepsilon}(W,h,r) + 2\|W\|_r(nV_h)^{\frac 1r- \frac 12}\right\}\\
     &\leq C_r^* \left\{\Lambda_{\varepsilon}(W,h,r) + 2\|W\|_r\right\}.
 \end{align}
 We have
 \begin{align*}\Lambda_\varepsilon(W,h,r)
 \le & F_\infty^{1/2} \left(\int_{\domain{d,\varepsilon}} \left(V_h\int_{\domain{d}} \mathcal K_{W,h}^2(t,x) \D x\right)^{r/2} \D t\right)^{1/r}\\
 \le & F_\infty^{1/2} \|W\|_2.
   \end{align*}

   As a consequence, for all $r\ge 1$, we have
   \begin{equation*}
   \Gamma_\varepsilon(W,h, r)\le C \|W\|_{r\vee 2}
   \end{equation*}
   where $C$ depends on $F_\infty$ and $r$.

 \subsection{Proof of Lemma~\ref{lem:bousq}}

Let $W\in\mathcal{W}^d$ and $h\in \Hn$.
 We denote by $\mathbb B_{r'}$ the unit ball of $\lebesgue_{r'}(\Delta_{d,\varepsilon})$ where $1/r+1/r'=1$ and, for $\lambda\in\mathbb B_{r'}$, we consider $\bar g_{\lambda}$ defined, for $x\in\Delta_d$ by:
 \begin{equation}
     \bar g_{\lambda}(x) = g_{\lambda}(x) - \e_f^n g_\lambda(X_1)
     \quad\text{with}\quad
     g_{\lambda}(x) = V_h^{1/2}\int_{\Delta_{d,\varepsilon}} \lambda(t) \mathcal K_{W,h}(t,x) \D t.
 \end{equation}
 The variable $Y= \|{\xi_{W,h}}\|_{r,\varepsilon}$ satisfies
 \begin{align}
    Y &= \sup_{\|\lambda\|_{r',\varepsilon}\leq 1} \int_{\Delta_{d,\varepsilon}} \lambda(t)\xi_{W,h}(t) \D t\\
     &= \sup_{\|\lambda\|_{r',\varepsilon}\leq 1} \frac1{\sqrt{n}}\sum_{j=1}^n \bar g_{\lambda}(X_j)
 \end{align}

 Since the set $\mathbb B_{r'}$ is a weakly--$*$ separable space, there exists a countable set {$(\lambda_k)_{k\in\N}\in \mathbb B_{r'}$} such that
 \begin{equation}
     Y = \sup_{k\in\N} \frac1{\sqrt{n}}\sum_{j=1}^n \bar g_{\lambda_k}(X_j).
 \end{equation}
{
We have
 \begin{align}
   \sup_{k\in\N} \|{\bar g_{\lambda_k}}\|_{\infty} &\leq 2\sup_{k\in\N}\|{g_{\lambda_k}}\|_{\infty} \\
    &\leq 2 \sup_{x\in\Delta_d} V_h^{1/2}\sup_{k\in\N}\|{\lambda_k}\|_{r',\varepsilon}\|{\mathcal K_{W,h}(\cdot,x)}\|_{r,\varepsilon}\\
    &\leq \mathfrak b(W,h,r),\label{eq:b}
 \end{align}
where
\begin{equation}
   \mathfrak b(W,h,r) =\mathfrak b=2\|W\|_r V_h^{1/r-1/2}.
 \end{equation}
For $r<2$, using the H\"older inequality, we have
 \begin{align}
     \sup_{k\in\N} \e_f^n g_{\lambda_k}^2(X_1)
     &= V_h \sup_{k\in\N} \int_{\Delta_d} \left(\int_{\Delta_{d,\varepsilon}} \lambda_k(t)\mathcal K_{W,h}(t,x) \D t\right) ^2 f(x) \D x\\
     &\leq V_h\sup_{k\in\N} \int_{\Delta_d}
     \|{\mathcal K_{W,h}(\cdot,x)}\|_{r,\varepsilon}^2
     \|{\lambda_k}\|_{r',\varepsilon}^2 f(x) \D x\\
     &= V_h^{2/r-1}
     \|W\|_r^2\label{eq:intermediaire1}.
 \end{align}
 For $r\geq 2$, using the Young inequality, we have
 \begin{align}
     \sup_{k\in\N} \e_f^n g_{\lambda_k}^2(X_1)
     &\leq
     \mathrm F_\infty V_h
     \sup_{k\in\N}
     \int_{\Delta_d}
     \left(
     \int_{\Delta_{d,\varepsilon}}
     \mathcal K_{W,h}(t,x)
     \lambda_k(t)
     \D t
     \right)^2 \D x\\
     &\leq \mathrm F_\infty V_h^{2/r}   \|W\|_{2r/(r+2)}^2\label{eq:intermediaire2}.
 \end{align}
 Combining \eqref{eq:intermediaire1} and \eqref{eq:intermediaire2}, we deduce that
 \begin{equation}\label{eq:sigma}
     \sup_{k\in\N} \e_f^n g_{\lambda_k}^2(X_1)\le \sigma^2(W,h,r),
 \end{equation}
 where
 \begin{equation}
   \sigma^2(W,h,r) =\sigma^2=
   \left\{\begin{array}{cc}
     \|W\|_r^2 V_h^{\frac2r-1} & \text{if $1 \leq r < 2$ }\\
     \mathrm F_\infty \|W\|_{2r/(r+2)}^2  V_h^{\frac2r}
     &\text{if $r\geq 2$}.
   \end{array}\right.
 \end{equation}}
 Now, using the Bousquet inequality \citep[see][]{bousquet2002bennett}, and denoting $\Gamma_{\varepsilon}=\Gamma_{\varepsilon}(W,h,r)$, we obtain for any $x>0$:
 \begin{align}
   \p\left(Y-\e_f^n Y \geq \frac{\Gamma_{\varepsilon}\tau}{2}+x\right)
   &\leq \exp\left(-\frac{x^2}{2\sigma^2 + \frac{\kb}{\sqrt{n}}\left(\Gamma_{\varepsilon}\left(\frac{12+\tau}{3}\right)+\frac{2x}{3}\right)}\right)\\
   &\qquad \times  \exp \left(-\frac{\tau\Gamma_{\varepsilon} x + \Gamma_{\varepsilon}^2\tau^2/4}{2\sigma^2 + \frac{\kb}{\sqrt{n}}\left(\Gamma_{\varepsilon}\left(\frac{12+\tau}{3}\right)+\frac{2x}{3}\right)}\right)\label{homo}
 \end{align}
Note that, for any $x>0$, we have
 \begin{equation}
 \frac{\tau\Gamma_{\varepsilon} x + \Gamma_{\varepsilon}^2\tau^2/4}{2\sigma^2 + \frac{\kb}{\sqrt{n}}\left(\Gamma_{\varepsilon}\left(\frac{12+\tau}{3}\right)+\frac{2x}{3}\right)}
 \geq
 \frac{\Gamma_{\varepsilon}^2\tau^2}{4\left(2\sigma^2 + \frac{\kb\Gamma_{\varepsilon}(12+\tau)}{3\sqrt{n}}\right)}.
 \end{equation}
This inequality holds due to the fact that the homography on the left hand side of the equation is an increasing function.
{Using \eqref{eq:b}, \eqref{eq:sigma} and the fact that if $r\le r'$, $\|W\|_r\le \|W\|_{r'}$, we obtain that for $r<2$
 \begin{align}
   \frac{4\left(2\sigma^2 + \frac{\kb\Gamma_{\varepsilon}(12+\tau)}{3\sqrt{n}}\right)}{\Gamma_{\varepsilon}^2\tau^2}
   &\leq \kc_1  (V_r)^{\frac{2}{r}-1} + \kc_2 (V_h)^{\frac{1}{r}-\frac{1}{2}}/\sqrt{n}\\
   &\leq \kc_1(h_n^*)^{d\left(\frac{2}{r}-1\right)} + \kc_2 (h_n^*)^{\frac{d}{r}}\label{eq:191}
 \end{align}
 where $\kc_1$ and $\kc_2$ are absolute positive constants that depend only on $d$, $\tau$ and $r$. For $r\ge 2$, using Lemma~\ref{lem:espY}, we have in a similar way
\begin{align}
   \frac{4\left(2\sigma^2 + \frac{\kb\Gamma_{\varepsilon}(12+\tau)}{3\sqrt{n}}\right)}{\Gamma_{\varepsilon}^2\tau^2}
   &\leq \kc_3  (V_r)^{\frac{2}{r}} + \kc_4 (V_h)^{\frac{1}{r}-\frac{1}{2}}/\sqrt{n}\\
   &\leq \kc_3(h_n^*)^{\frac{2d}{r}} + \kc_4 (h_n^*)^{\frac{d}{r}}.\label{eq:192}
 \end{align}
  where $\kc_3$ and $\kc_4$ are absolute positive constants that depend only on $d$, $\tau$, $F_\infty$ and $r$.
 Finally, we deduce that
 \begin{align}
   \frac{4\left(2\sigma^2 + \frac{\kb\Gamma_{\varepsilon}(12+\tau)}{3\sqrt{n}}\right)}{\Gamma_{\varepsilon}^2\tau^2}\label{eq:19}
   &\leq \kc_5  \left(\alpha_n(r)\right)^{-1},
  \end{align}
 with $\kc_5$ an absolute positive constant that depends only on $d$, $\tau$, $F_\infty$ and $r$.
 Using \eqref{eq:19} we obtain:
 \begin{equation}
   \label{eq:18}
   \exp \left(-\frac{\tau\Gamma_{\varepsilon} x + \Gamma_{\varepsilon}^2\tau^2/4}{2\sigma^2 + \frac{\kb}{\sqrt{n}}\left(\Gamma_{\varepsilon}\left(\frac{12+\tau}{3}\right)+\frac{2x}{3}\right)}\right)
 \leq \exp(-C_1\alpha_n(r)),
 \end{equation}
 where $C_1$ is an absolute positive constant that depends only on $r$, $\tau$ and $\mathrm F_\infty$.
 \ \\
Using Lemma \ref{lem:espY}, \eqref{eq:b} and \eqref{eq:sigma}, we obtain that there exists an absolute constant $\kc_6$ which depends only on $\mathrm F_\infty$, $\tau$ and $r$ such that:
 \begin{equation}
   \label{eq:4}
   2\sigma^2 + \frac{\kb}{\sqrt{n}}\left(\Gamma_{\varepsilon}\left(\frac{12+\tau}{3}\right)+\frac{2x}{3}\right)
   \leq \kc_6\alpha_n(r)(\|W\|^2_{2\vee r}+x\|W\|_{r})
 \end{equation}
\eqref{homo}, \eqref{eq:18} and \eqref{eq:4}, allow us to deduce the result of the lemma.}

\subsection{Proof of Lemma~\ref{lem:bias1}.}
Let $(e_1,\ldots,e_d)$ be the canonical basis of $\R^d$ and define
\begin{equation}
    v_i(u) = (t_1+\eta_1u_1,\ldots,t_{i-1}+\eta_{i-1}u_{i-1},
    \quad t_i\quad,
    t_{i+1}+h_{i+1}u_{i+1}, \ldots, t_d+h_du_d).
\end{equation}
We can write:
\begin{align}
    f(t+h\cdot u)-f(t+\eta\cdot u)
    &= \sum_{i=1}^d f(v_i(u)+h_iu_ie_i) - f(v_i(u)+\eta_iu_ie_i)\\
    &= \sum_{i\in I} f(v_i(u)+h_iu_ie_i) - f(v_i(u)),
\end{align}
where $I=\{i=1,\ldots,d : \eta_i=0\}$. Using a Taylor expansion of the function $x\in\R\mapsto f(v_i(u)+xe_i)$ around $0$, we obtain:
\begin{align}
    f(t+h\cdot u)-&f(t+\eta\cdot u)
    = \sum_{i\in I}
        {\sum_{k=1}^{\lfloor s_i\rfloor}}
        D_i^k f(v_i(u))\frac{(h_iu_i)^k}{k!}\\
    &+\sum_{i\in I} \frac{(h_iu_i)^{\lfloor s_i\rfloor}}{\lfloor s_i\rfloor!}
        \int_0^1 {(1-\tau)^{\lfloor s_i\rfloor-1}}\bigl[D_i^{\lfloor s_i\rfloor} f(v_i(u)+\tau h_iu_i)-D_i^{\lfloor s_i\rfloor} f(v_i(u))\bigr] \D \tau.
\end{align}
Using the facts that $v_i(u)$ does not depend on $u_i$ and that $\int_{\Delta_1} W_i(y) y^k \D y=0$ for any $1\leq k\leq\lfloor s_i\rfloor$, Fubini's theorem implies that:
\begin{align}
    S_{W,h, \eta}^*(f)
    &= \left(
    \int_{\Delta_{d,\mathbf 0}}
    \left\vert
    \int_{\Delta_d}
    \left(\prod_{i=1}^d W_i(u_i)\right)
    \sum_{i\in I}
    {I_i(t,u,h)}
    \D u
    \right\vert^p
    \D t
    \right)^{1/p}
\end{align}
where
\begin{equation}
    I_i(t,u,h) = \frac{(h_iu_i)^{\lfloor s_i\rfloor}}{\lfloor s_i\rfloor!}
    \int_0^1 {(1-\tau)^{\lfloor s_i\rfloor-1}}\bigr[D_i^{\lfloor s_i\rfloor} f(v_i(u)+\tau h_iu_i)-D_i^{\lfloor s_i\rfloor} f(v_i(u))\bigl] \D \tau.
\end{equation}
Using Jensen's inequality and Fubini's theorem we obtain that:
\begin{align}
    S_{W,h, \eta}^*(f)
    &=
    \big(d\|W\|_1\big)^{1-1/p}
    \left(
    \int_{\Delta_d}
    J(u,h)
    \biggl\vert
    \prod_{i=1}^d W_i(u_i)
    \biggr\vert
    \D u
    \right)^{1/p},
\end{align}
where $J(u,h) = \sum_{i\in I}
    \int_{\Delta_{d,\mathbf 0}}
    \left\vert
    I_i(t,u,h)
    \right\vert^p
    \D t$. Now, we study this last term:
\begin{align}
    J(u,h) &\leq
    \sum_{i\in I}
    \int_{\Delta_{d,\mathbf 0}}
    \frac{(h_iu_i)^{1+p s_i}}{(\lfloor s_i\rfloor!)^p}
    \int_0^1
    \frac{\left|D_i^{\lfloor s_i\rfloor} f(v_i(u)+\tau h_iu_i)-D_i^{\lfloor s_i\rfloor} f(v_i(u))\right|^p}%
    {|\tau h_iu_i|^{1+p(s_i-\lfloor s_i\rfloor)}}
    \D\tau
    \D t.
\end{align}
Using a simple change of variables, we obtain:
\begin{align}
    J(u,h) &\leq
    \sum_{i\in I}
    \frac{(h_iu_i)^{p s_i}}{(\lfloor s_i\rfloor!)^p}
    \int_{\Delta_{d}}
    \int_0^1
    \frac{\left|D_i^{\lfloor s_i\rfloor} f(x_1,\ldots,x_{i-1},\xi,x_{i+1},\ldots,x_d)-D_i^{\lfloor s_i\rfloor} f(x)\right|^p}%
    {|\xi-x_i|^{1+p(s_i-\lfloor s_i\rfloor)}}
    \D\xi
    \D x.
\end{align}
Since $u_i\leq 1$ and $f\in\sobolev_{s,p}(L)$ we have:
\begin{equation}
    S_{W,h, \eta}^*(f) \leq d\|W\|_1
    \kappa(s)L\left(
    \sum_{i\in I}  h_i^{ps_i}
    \right)^{1/p}.
\end{equation}
This implies the result.



\section{Acknowledgments}
The authors have been supported by Fondecyt project 1171335, Mathamsud project 18-MATH-07 and ECOS project C15E05.
\bibliographystyle{plainnat}
\bibliography{statnew}
\end{document}